\documentclass[12pt]{article}
\usepackage{amstex,amssymb,amscd}
\begin{document}
%%%%%%%%%%%%%%%%%%%  Environments  %%%%%%%%%%%%%%%%%%%%%%%%%%%%%
\newtheorem{thm}{Theorem}[section]
\newtheorem{lem}[thm]{Lemma}
\newtheorem{prop}[thm]{Proposition}
\newtheorem{conj}[thm]{Conjecture}
\newtheorem{cor}[thm]{Corollary}
\newenvironment{dfn}{\medskip\refstepcounter{thm}
\noindent{\bf Definition \thesection.\arabic{thm}\ }}{\medskip}
\newenvironment{ex}{\medskip\refstepcounter{thm}
\noindent{\bf Example \thesection.\arabic{thm}\ }}{\medskip}
\newenvironment{proof}{\medskip\noindent{\it Proof.\ }}
{\hfill\hskip 2em$\Box$\par}
%%%%%%%%%%%%%%%%%%%%%%%  Macros  %%%%%%%%%%%%%%%%%%%%%%%%%%%%%%%
\def\dim{\mathop{\rm dim}}
\def\Hol{\mathop{\rm Hol}}
\def\Ric{\mathop{\rm Ric}}
\def\id{\mathop{\rm id}}
\def\Fix{\mathop{\rm Fix}}
\def\eq#1{{\rm(\ref{#1})}}
\def\R{\mathbin{\mathbb R}}
\def\C{{\textstyle\mathop{\mathbb C}}}
\def\d{{\rm d}}
\def\ms#1{\vert#1\vert^2}
\def\bms#1{\bigl\vert#1\bigr\vert^2}
\def\nm#1{\Vert #1 \Vert}
\def\cnm#1#2{\Vert #1 \Vert_{C^{#2}}}
\def\md#1{\vert #1 \vert}
\def\bnm#1{\bigl\Vert #1 \bigr\Vert}
\def\bmd#1{\bigl\vert #1 \bigr\vert}
%%%%%%%%%%%%%%%%%%  Text of paper  %%%%%%%%%%%%%%%%%%%%%%%%%%%%%
\title{Quasi-ALE metrics with \\ holonomy ${\rm SU}(m)$ and ${\rm Sp}(m)$}
\author{Dominic Joyce \\
Lincoln College, Oxford, OX1 3DR, England}
\date{May 1999}
\maketitle

This is the sequel to the author's paper \cite{Joyc5} on 
{\it Asymptotically Locally Euclidean metrics}, or {\it ALE metrics} 
for short, with holonomy ${\rm SU}(m)$. Let $G$ be a finite subgroup of 
${\rm U}(m)$ and suppose $(X,\pi)$ is a {\it resolution} of $\C^m/G$, 
that is, $X$ is a normal nonsingular variety with a proper birational 
morphism~$\pi:X\rightarrow\C^m/G$.

When $G$ acts freely on $\C^m\setminus\{0\}$, so that $\C^m/G$ has 
an {\it isolated quotient singularity} at 0, we defined in \cite{Joyc5} 
a special class of K\"ahler metrics on $X$ called {\it ALE K\"ahler 
metrics}, and proved an existence result for Ricci-flat ALE K\"ahler 
metrics on {\it crepant resolutions} of $\C^m/G$, using a version of 
the Calabi conjecture for ALE manifolds.

This paper will generalize the ideas of \cite{Joyc5} to 
the case when $G$ does not act freely on $\C^m\setminus\{0\}$, 
so that the singularities of $\C^m/G$ are not isolated. 
The appropriate class of K\"ahler metrics on resolutions $X$ of 
non-isolated quotient singularities $\C^m/G$ will be called 
{\it Quasi-ALE}, or {\it QALE} for short. We are particularly
interested in {\it Ricci-flat} QALE K\"ahler manifolds, and our
main result is an existence theorem for Ricci-flat QALE K\"ahler 
metrics on crepant resolutions of~$\C^m/G$.

The key to understanding the structure of non-isolated singularities
$\C^m/G$ is to observe that if $s$ is a singular point of $\C^m/G$
and $s\ne 0$, then an open neighbourhood of $s$ in $\C^m/G$ is
isomorphic to an open neighbourhood of $(0,0)$ in $\C^k\times\C^{m-k}/H$,
where $0<k<m$ and $H$ is a finite subgroup of $U(m-k)$, and also of 
$G$. Thus, away from zero the singularities of $\C^m/G$ look locally 
like {\it products} $\C^k\times\C^{m-k}/H$ for~$k>0$. 

In \S 1 we will define a special class of resolutions $X$ 
of $\C^m/G$ called {\it local product resolutions}, which have the 
property that if $\C^m/G$ is locally modelled on $\C^k\times\C^{m-k}/H$ 
then $X$ is locally modelled on $\C^k\times Y$, where $Y$ is a resolution 
of $\C^{m-k}/H$. Crepant resolutions are automatically of this form.
If $X$ is a local product resolution, we want to impose some suitable 
asymptotic conditions `near infinity' on K\"ahler metrics $g$ on~$X$. 

In \S 2 we define $g_{\scriptscriptstyle X}$ to be a 
{\it QALE K\"ahler metric} on $X$ if $g_{\scriptscriptstyle X}$ 
is asymptotic to $h_{\smash{\scriptscriptstyle{\mathbb C}^k}}\times 
g_{\scriptscriptstyle Y}$ on the part of $X$ modelled on $\C^k\times Y$, 
where $h_{\smash{\scriptscriptstyle{\mathbb C}^k}}$ is the Euclidean 
metric on $\C^k$, and $g_{\scriptscriptstyle Y}$ is a metric on $Y$, 
the resolution of $\C^{m-k}/H$. So $g_{\scriptscriptstyle X}$ 
converges to $h_{\smash{\scriptscriptstyle{\mathbb C}^k}}\times 
g_{\scriptscriptstyle Y}$ at infinity. Now the points of $\C^m/G$ 
will be of several different kinds, modelled on $\C^k\times\C^{m-k}/H$ 
for different $0\le k\le m$ and subgroups $H\subseteq G$. Thus we impose 
not one but many asymptotic conditions on $g_{\scriptscriptstyle X}$, 
which must all be satisfied for $g_{\scriptscriptstyle X}$ to be QALE.

Section 3 discusses {\it Ricci-flat QALE K\"ahler manifolds}.
We state the main result of the paper, Theorem \ref{qalerfthm}, 
an existence result for Ricci-flat QALE metrics on crepant 
resolutions $X$ of $\C^m/G$. Its proof follows that of Theorem 
3.3 of \cite{Joyc5}, and takes up most of sections 4-7. 
We first develop the appropriate ideas of weighted H\"older 
spaces and elliptic regularity on QALE manifolds. Then 
we state (without proof) two versions of the Calabi conjecture 
on QALE manifolds, and apply them to construct Ricci-flat QALE 
K\"ahler metrics on~$X$.

The original motivation for this paper and \cite{Joyc5} is that
ALE and QALE metrics with holonomy SU(2), SU(3), SU(4) and Sp(2)
are essential ingredients in a new construction by the author
of compact manifolds with the exceptional holonomy groups $G_2$
and Spin(7), which generalizes that of \cite{Joyc1,Joyc2}. This
construction will be described at length in the author's 
forthcoming book~\cite{Joyc4}. 

For simplicity we restrict our attention in this paper to
{\it resolutions} of $\C^m/G$. But in \cite[\S 9.9]{Joyc4} we 
will give a more general definition of QALE K\"ahler manifolds, 
in which the underlying complex manifold is allowed to be a 
{\it deformation} of $\C^m/G$, or a resolution of a deformation. 
We show that our main result, Theorem \ref{qalerfthm}, also 
holds in this more general context.

The book will also contain the results of this paper and \cite{Joyc5}, 
and many other results on QALE manifolds, including the proofs of 
the QALE Calabi conjectures stated in \S 6, and existence results 
for QALE manifolds with holonomy $G_2$ and Spin(7). The author 
does not know of any previous papers on QALE manifolds at all,
at the time of writing; this may be the first paper on the subject.

\section{Local product resolutions}

In this section we will study the structure of {\it nonisolated 
quotient singularities} $\C^m/G$ and their resolutions. Here a {\it 
resolution} $(X,\pi)$ of $\C^m/G$ is a normal nonsingular variety 
$X$ with a proper birational morphism $\pi:X\rightarrow\C^m/G$. We 
shall define a special kind of resolution of $\C^m/G$ called a
{\it local product resolution}. Our main goal is to set up a lot 
of notation describing resolutions of $\C^m/G$, that we will use 
in the rest of the paper.

Let $G$ be a finite subgroup of ${\rm U}(m)$. If $A$ is a subgroup of 
$G$ and $V$ a subspace of $\C^m$, define the {\it fixed point set}\/ 
$\Fix(A)$, the {\it centralizer} $C(V)$ and the {\it normalizer} $N(V)$ by
\begin{equation}
\begin{split}
\Fix(A)&=\{x\in\C^m:\text{$a\,x=x$ for all $a\in A$}\},\\
C(V)&=\{g\in G:\text{$g\,v=v$ for all $v\in V$}\}\\
\text{and}\quad N(V)&=\{g\in G:g\,V=V\}.
\end{split}
\end{equation}
Then $C(V)$ and $N(V)$ are subgroups of $G$, and $C(V)$
is a normal subgroup of~$N(V)$.

\begin{dfn} Define a finite set $\mathcal L$ of linear subspaces of 
$\C^m$ by
\begin{equation}
{\mathcal L}=\bigl\{\Fix(A):\text{$A$ is a subgroup of $G$}\bigr\}.
\end{equation}
Let $I$ be an {\it indexing set}\/ for $\mathcal L$, so that we 
may write ${\mathcal L}=\{V_i:i\in I\}$. Let the indices of 
$\Fix\bigl(\{1\}\bigr)$ and $\Fix(G)$ be $0,\infty$ respectively, 
so that $0,\infty\in I$ and $V_0=\C^m$, $V_\infty=\Fix(G)$ by 
definition. Usually $V_\infty=\{0\}$. 

Define a {\it partial order} $\succeq$ on $I$ by $i\succeq j$ if 
$V_i\subseteq V_j$. Then $\infty\succeq i\succeq 0$ for all $i\in I$. 
Let $W_i$ be the perpendicular subspace to $V_i$ in $\C^m$, so 
that $\C^m=V_i\oplus W_i$. Define $A_i=C(V_i)$. Then $V_i=\Fix(A_i)$ 
and $A_i$ acts on $W_i$, with $\C^m/A_i\cong V_i\times W_i/A_i$. If 
$i\succeq j$ then $W_i\supseteq W_j$ and $A_i\supseteq A_j$. 
Define $B_i$ to be the quotient group $N(V_i)/C(V_i)$. Then 
$B_i$ acts naturally on $V_i$ and $W_i/A_i$, and 
$(V_i\times W_i/A_i)/B_i\cong\C^m/N(V_i)$. Hence, if $N(V_i)=G$ 
then~$(V_i\times W_i/A_i)/B_i\cong\C^m/G$.
\label{liviwidef}
\end{dfn}

Let $V_i,V_j\in\mathcal L$, let $A$ be the subgroup of $G$ 
generated by $A_i$ and $A_j$, and set $V=\Fix(A)$. It is easy 
to show that $V=V_i\cap V_j$. Thus $V_i\cap V_j\in\mathcal L$, and 
$\mathcal L$ is closed under intersection of subspaces. Also, if
$g\in G$ and $V_i\in\mathcal L$ then $g\,V_i\in\mathcal L$, since 
$g\,V_i=\Fix(g\,A_ig^{-1})$. For each $g\in G$ and $i\in I$, let 
$g\cdot i$ be the unique element of $I$ such that $V_{g\cdot i}
=g\,V_i$. This defines an action of $G$ on $I$, which satisfies 
$W_{g\cdot i}=g\,W_i$ and~$A_{g\cdot i}=g\,A_ig^{-1}$.

Let $v\in\C^m$. Then $vG$ is a singular point of $\C^m/G$ if and 
only if the subgroup of $G$ fixing $v$ is nontrivial, that is, if 
$v\in\Fix(A)$ for some nontrivial subgroup $A\subset G$. Thus $vG$ 
is a singular point if and only if $v\in V_i$ for some $i\in I$ with 
$i\ne 0$, and the {\it singular set} $S$ of $\C^m/G$ is
\begin{equation}
S=\bigcup\begin{Sb}i\in I\setminus\{0\}\end{Sb}V_i\Bigl/G.
\end{equation}
For generic points $v\in V_i$, the subgroup of $G$ fixing $v$ is 
$A_i$, and the singularity of $\C^m/G$ at $vG$ is locally modelled 
on the product~$V_i\times W_i/A_i$. 

\begin{dfn} Let $(X,\pi)$ be a resolution of $\C^m/G$. We say 
that $(X,\pi)$ is a {\it local product resolution} if for each 
$i\in I$ there exists a resolution $(Y_i,\pi_i)$ of $W_i/A_i$ 
such that the following conditions hold. Let $R>0$, and
define subsets $S_i$ and $T_i$ in $V_i\times W_i/A_i\cong\C^m/A_i$ by
\begin{align}
S_i&=\bigcup\begin{Sb}j\in I:i\nsucceq j\end{Sb}V_j\,\biggl/A_i,\\
\begin{split}
T_i&=\bigl\{x\in\C^m/A_i:d(x,S_i)\le R\bigr\}\\
&=\bigl\{x\in\C^m:\text{$d(x,V_j)\le R$ for some $j\in I$
with $i\nsucceq j$}\bigr\}\Big/A_i,
\end{split}
\end{align}
where $d(\,,\,)$ is the distance in $\C^m/A_i$ or $\C^m$. 
Let $U_i$ be the pull-back of $T_i$ to $V_i\times Y_i$ under 
$\id\times\pi_i:V_i\times Y_i\rightarrow V_i\times W_i/A_i$. 
Let $\phi_i$ be the natural projection from $V_i\times W_i/A_i$ 
to $\C^m/G$. Then there should exist a map $\psi_i:V_i\times 
Y_i\big\backslash U_i\rightarrow X$ such that the following 
diagram commutes:
\begin{equation}
\begin{CD}
V_i\times Y_i\,\Big\backslash U_i@>{\psi_i}>> X\\
@VV{\id\times\pi_i}V @VV{\pi}V\\
V_i\times W_i/A_i\,\Big\backslash T_i @>{\phi_i}>>
{\mathbin{\mathbb C}}^m/G.
\end{CD}
\label{phipicdeq}
\end{equation}
This $\psi_i$ should be a local isomorphism of complex manifolds,
and whenever $x\in X$ and $y\in V_i\times W_i/A_i\big\backslash T_i$
satisfy $\pi(x)=\phi_i(y)$, there should exist a unique
$z\in V_i\times Y_i\big\backslash U_i$ with $x=\psi_i(z)$ 
and~$y=(\id\times\pi_i)(z)$. 
\label{locproddef}
\end{dfn}

Here is an equivalent way to write this. Define $\tilde X_i$ by
\begin{equation*}
\tilde X_i=\Bigl\{(x,y)\in X\times\bigl(V_i\times W_i/A_i\big\backslash 
T_i\bigr):\pi(x)=\phi_i(y)\Bigr\}.
\end{equation*}
For any resolution $X$ of $\C^m/G$, one can show that $\tilde X_i$ is 
a well-defined complex manifold, and the projection to $X$ is a local 
isomorphism. We say that $X$ is a local product resolution if $\tilde X_i$ 
is isomorphic to $V_i\times Y_i\big\backslash U_i$ in such a way that 
the natural projections to $V_i\times W_i/A_i\big\backslash T_i$ agree.

The idea of the definition is that every point in $\C^m/G$ has an 
open neighbourhood isomorphic to an open neighbourhood of 0 in 
$V_i\times W_i/A_i$ for some $i\in I$. (For nonsingular points we 
use $V_0\times W_0/A_0=\C^m$.) A resolution $X$ of $\C^m/G$ is a 
local product resolution if and only if the resolution of each 
singular point modelled on $V_i\times W_i/A_i$ is locally isomorphic 
to~$V_i\times Y_i$.

Local product resolutions are a very large class of resolutions
of $\C^m/G$. There do exist (rather artificial) examples of
resolutions which are not local product resolutions, but in
practice they include all interesting resolutions of $\C^m/G$,
in particular crepant resolutions, as we shall see in \S 3.

\begin{prop} Let\/ $(X,\pi)$ be a local product resolution
of\/ $\C^m/G$. Then for each\/ $g\in G$ and\/ $i\in I$ there 
is a unique isomorphism $\chi_{g,i}:Y_i\rightarrow Y_{g\cdot i}$
making a commutative diagram
\begin{equation}
\begin{CD}
Y_i@>{\chi_{g,i}}>> Y_{g\cdot i}\\
@VV{\pi_i}V @VV{\pi_{g\cdot i}}V\\
W_i/A_i@>{g}>>W_{g\cdot i}/A_{g\cdot i},
\end{CD}
\label{chigicdeq}
\end{equation}
where the map $W_i/A_i{\buildrel g\over\longrightarrow}\,W_{g\cdot i}/
A_{g\cdot i}$ is given by~$g(w\,A_i)=(g\,w)A_{g\cdot i}$.
\label{chigiprop}
\end{prop}

\begin{proof} If $v$ is a generic point of $V_i$ then $\C^m/G$ is 
locally isomorphic to $V_i\times W_i/A_i$ near $vG$, and $X$ is locally 
isomorphic to $V_i\times Y_i$ near $\pi^{-1}(vG)$, as $X$ is a local 
product resolution. But if $g\in G$ then $vG=(gv)G$, and $gv$ is a 
generic point of $V_{g\cdot i}$, so that $\C^m/G$ is also locally 
isomorphic to $V_{g\cdot i}\times W_{g\cdot i}/A_{g\cdot i}$ near $vG$, 
and $X$ is locally isomorphic to $V_{g\cdot i}\times Y_{g\cdot i}$ near
$\pi^{-1}(vG)$. Hence $V_i\times Y_i$ and $V_{g\cdot i}\times Y_{g\cdot i}$
are locally isomorphic. In fact they are globally isomorphic, and 
this gives an isomorphism $\chi_{g,i}:Y_i\rightarrow Y_{g\cdot i}$, which 
makes \eq{chigicdeq} commutative.
\end{proof}

If $g\in N(V_i)$ then $g\cdot i=i$, so that $\chi_{g,i}$ is an
automorphism of $Y_i$. Also, if $g\in A_i$ then $\chi_{g,i}$ is
the identity on $Y_i$. Thus $N(V_i)$ acts on $Y_i$, and the
normal subgroup $A_i$ of $N(V_i)$ acts trivially on $Y_i$.
Therefore the action descends to an action of the quotient
group $B_i=N(V_i)/A_i$ on $Y_i$. That is, the action of
$B_i$ on $W_i/A_i$ must lift to an action of $B_i$ on the
resolution $Y_i$ on $W_i/A_i$. So each resolution $Y_i$ in 
a local product resolution must be $B_i$-equivariant.

We can use this to give a more thorough explanation of 
the idea of local product resolution, under a simplifying 
assumption. Suppose that $N(V_i)=G$ for each $i\in I$. This 
happens if $G$ is abelian, and in other cases too. Then we 
have $(V_i\times W_i/A_i)/B_i\cong\C^m/G$ for $i\in I$, from
Definition \ref{liviwidef}. Since $B_i$ acts on $V_i$ and
$Y_i$ we can take the quotient $(V_i\times Y_i)/B_i$, which is 
a complex orbifold with a natural projection to~$(V_i\times 
W_i/A_i)/B_i\cong\C^m/G$. 

Effectively, $(V_i\times Y_i)/B_i$ is a {\it partial resolution} 
of $\C^m/G$, which resolves the singularities of $\C^m/G$ due 
to fixed points of elements of $A_i$, but leaves unresolved 
those caused by fixed points of elements of $G\setminus A_i$. Now
$S_i$ is exactly the subset of $V_i\times W_i/A_i$ fixed by some 
$b\ne 1$ in $B_i$, and therefore $S_i/B_i$ is the set of 
singularities of $\C^m/G$ due to fixed points of $G\setminus A_i$. 
Also $T_i$ is the subset of $V_i\times W_i/A_i$ within distance 
$R$ of $S_i$, so $T_i/B_i$ is the subset of $\C^m/G$ within 
distance $R$ of singularities due to fixed points of~$G\setminus A_i$.

All the fixed points of the $B_i$-action on $V_i\times Y_i$ 
lie in $(\id\times\pi_i)^{-1}(S_i)$, which is in the interior 
of $U_i$. Thus $B_i$ acts freely on $V_i\times Y_i\setminus U_i$, 
and $(V_i\times Y_i\setminus U_i)/B_i$ is nonsingular. In fact 
$\psi_i:(V_i\times Y_i\setminus U_i)/B_i\rightarrow X$ is an 
{\it isomorphism} between $(V_i\times Y_i\setminus U_i)/B_i$ 
and $X\setminus\pi^{-1}(T_i/B_i)$. Thus, the definition requires 
that the resolution $X$ of $\C^m/G$ has to coincide with the 
partial resolution $(V_i\times Y_i)/B_i$ of $\C^m/G$ outside the set 
$U_i/B_i$ which contains the singularities of~$(V_i\times Y_i)/B_i$.

Our next result shows that local product resolutions are built 
out of local product resolutions of smaller dimension.

\begin{prop} Suppose $X$ is a local product resolution of\/ 
$\C^m/G$. Then each of the resolutions $Y_i$ of\/ $W_i/A_i$ in 
Definition \ref{locproddef} is also a local product resolution.
\label{lpprop}
\end{prop}

\begin{proof} Let $i\in I$ be fixed throughout the proof. Above 
we defined a lot of notation such as ${\mathcal L},I,V_j,W_j$ and 
so on, associated to $\C^m/G$. The corresponding data associated 
to $W_i/A_i$ will be written ${\mathcal L}',I',V_j',W_j'$, etc., 
in the obvious way. We will express the data for $W_i/A_i$ in 
terms of that for $\C^m/G$. Define the indexing set $I'$ by 
$I'=\{j\in I:i\succeq j\}$, and for each $j\in I'$ set
$V_j'=V_j\cap W_i$. Then ${\mathcal L}'=\{V_j':j\in I'\}=
\bigl\{\Fix(A):A$ is a subgroup of $A_i\bigr\}$ is a finite 
set of subspaces of $W_i$. The two special elements of $I'$ 
are $0'=0$ and $\infty'=i$. Also $W_j'$, $A_j'$ and $Y_j'$
are the same as $W_j$, $A_j$ and $Y_j$ for each~$j\in I'$.

Let $\phi_j'$, $T_j'$ and $U_j'$ be as in Definition 
\ref{locproddef}. It can be shown that there is a unique map
$\psi_j':V_j'\times Y_j'\big\backslash U_j'\rightarrow Y_i$, such 
that the product $\id\times\psi_j'$ of $\psi_j'$ with the identity 
on $V_i$ makes the following picture into a commutative diagram:
\begin{equation}
\begin{CD}
V_i\times V_j'\times Y_j'\,\Big\backslash U_j@>{\id\times\psi_j'}>>
V_i\times Y_i\,\Big\backslash U_i\\
@VV{\cong}V   @VV{\psi_i}V\\
V_j\times Y_j\,\Big\backslash U_j @>{\psi_j}>>X.
\end{CD}
\label{ijjpcdeq}
\end{equation}
It easily follows that $Y_i$ is a local product resolution 
of~$W_i/A_i$.
\end{proof}

\section{Quasi-ALE K\"ahler metrics}

We will now define a class of K\"ahler metrics on local product
resolutions of $\C^m/G$ called {\it Quasi Asymptotically 
Locally Euclidean}, or {\it Quasi-ALE} or {\it QALE} for 
short, which generalize the ALE K\"ahler metrics considered in 
\cite{Joyc5}. Let $G$ be a finite subgroup of ${\rm U}(m)$, 
let $X$ be a local product resolution of $\C^m/G$, and let 
all notation be as in~\S 1.

\begin{dfn} For each pair $i,j\in I$, define 
$\mu_{i,j}:V_i\times Y_i\rightarrow[0,\infty)$ by 
$\mu_{i,j}(z)=d\bigl((\id\times\pi_i)(z),V_jA_i/A_i\bigr)$, 
where $V_jA_i/A_i=\{vA_i:v\in V_j\}$, as a subset of $\C^m/A_i$, 
and $d(y,T)$ is the shortest distance between the point $y$ 
and the subset $T$ in $\C^m/A_i$. For each $i\in I$, define
$\nu_i:V_i\times Y_i\rightarrow[1,\infty)$ by 
$\nu_i(z)=1+\min\bigl\{\mu_{i,j}(z):j\in I$, $j\ne 0\bigr\}$. 
Then $\mu_{i,j}$ and $\nu_i$ are both continuous functions 
on $V_i\times Y_i$. For each $i\in I$ define 
$d_i=2-2\dim W_i=2-2m+2\dim V_i$, and let $h_i$ be the 
Euclidean metric on~$V_i$. 

Let $g$ be a K\"ahler metric on $X$. We say $g$ is {\it Quasi-ALE} 
if the complex codimension $n$ of the singular set $S$ of $\C^m/G$, 
given by $n=\min\{\dim W_i:i\in I$, $i\ne 0\}$ satisfies $n\ge 2$,
and for each $i\in I$ there is a K\"ahler metric $g_i$ on $Y_i$, 
such that the metric $h_i\times g_i$ on $V_i\times Y_i$ satisfies
\begin{equation}
\nabla^l\bigl(\psi_i^*(g)-h_i\times g_i\bigr)=
\sum\begin{Sb}j\in I:i\nsucceq j\end{Sb}
O\bigl(\mu_{i,j}^{d_j}\nu_i^{-2-l}\bigr)
\label{qaleeq}
\end{equation}
on $V_i\times Y_i\,\big\backslash U_i$, for all $l\ge 0$. If $X$ 
is a local product resolution of $\C^m/G$ with complex structure 
$J$, and $g$ is a QALE metric on $X$, then we say that $(X,J,g)$ 
is a {\it QALE K\"ahler manifold asymptotic to}~$\C^m/G$.
\label{qaledef}
\end{dfn}

We now discuss this definition. The pull-back to $V_i\times Y_i$
of the singular set $S$ of $\C^m/G$ splits up into a number
of pieces parametrized by $j\in I\setminus\{0\}$, and $\mu_{i,j}$ 
is a measure of the distance to piece $j$. Similarly, $\nu_i$ 
measures the distance in $V_i\times Y_i$ to the pull-back of $S$, 
but we add 1 to this to avoid problems when $\nu_i$ is small. 
Also, by definition $U_i$ is the subset of $V_i\times Y_i$ on 
which $\mu_{i,j}\le R$ for some $j\in I$ with $i\nsucceq j$. 
Thus $\mu_{i,j}>R>0$ in \eq{qaleeq}, so there are no problems 
when $\mu_{i,j}$ is small.

Thus eqn \eq{qaleeq} says that at large distances from $U_i$, 
the pull-back $\psi_i^*(g)$ of $g$ to $V_i\times Y_i$ must 
approximate the product metric $h_i\times g_i$ on $V_i\times Y_i$. 
As in \S 1, we can explain this more clearly if we assume 
that $N(V_i)=G$. In this case $h_i\times g_i$ pushes down to an 
orbifold metric on $(V_i\times Y_i)/B_i$, and $X$ is isomorphic to 
$(V_i\times Y_i)/B_i$ outside $U_i/B_i$. So \eq{qaleeq} says that 
the metrics $g$ and $h_i\times g_i$ on the isomorphic subsets of 
$X$ and $(V_i\times Y_i)/B_i$ must agree asymptotically at large 
distances from~$U_i/B_i$. 

We assume that the singular set $S$ of $\C^m/G$ has 
codimension $n\ge 2$ because many of the results in the rest of 
the paper are false when $n=1$. One reason for this is that if 
$\dim W_j=1$ then $d_j=2-2\dim W_j=0$, and so the `error term' 
$O(\mu_{i,j}^{d_j}\nu_i^{-2-l})$ in \eq{qaleeq} may not be small 
when $\mu_{i,j}$ is large. But some of our results rely on the 
errors being small at large distances. 

Another reason is that the equation $\Delta u=f$ on $\C^k$ behaves 
differently in the cases $k=1$ and $k>1$. When $k>1$ and $f$ is a 
smooth function on $\C^k$ that decays rapidly at infinity, there 
is a unique smooth function $u$ on $\C^k$ with $\Delta u=f$ and 
$u=O(r^{2-2k})$ for large $r$. But when $k=1$ this is false, and 
instead $u=O(\log r)$ and is not unique. Because of this, some of
our results about the Laplacian on QALE manifolds are false 
when~$n=1$. 

Consider what \eq{qaleeq} means when $i=0$ and $i=\infty$. Now 
$V_0=\C^m$ and $W_0$ and $Y_0$ are both a single point, so 
$V_0\times W_0\cong\C^m$. The metric $h_0\times g_0$ on 
$V_0\times W_0$ is just the Euclidean metric $h_0$ on $\C^m$, and 
\eq{qaleeq} says that the metrics $g$ on $X$ and $h_0$ on $\C^m/G$ 
must be asymptotic at large distances from the singular set of~$\C^m/G$.

When $i=\infty$ we have $U_\infty=\emptyset$, the map 
$\psi_\infty:V_\infty\times Y_\infty\rightarrow X$ is an 
isomorphism, and the r.h.s.\ of \eq{qaleeq} is zero, which forces 
$\psi_\infty^*(g)=h_\infty\times g_\infty$. Using $\psi_\infty$
to identify $X$ and $V_\infty\times Y_\infty$, we see that 
$g=h_\infty\times g_\infty$, the product of a Euclidean metric 
on $V_\infty$ and a metric $g_\infty$ on $Y_\infty$. If 
$\Fix(G)=\{0\}$ then $V_\infty=\{0\}$ and $Y_\infty=X$, and 
\eq{qaleeq} holds trivially by taking $g_\infty=g$. It is often
convenient to assume that~$\Fix(G)=\{0\}$.

In particular, suppose $\C^m/G$ has an {\it isolated} 
singularity at 0. Then $I=\{0,\infty\}$ and $\Fix(G)=\{0\}$,
so \eq{qaleeq} is trivial for $i=\infty$. We have $d_\infty=2-2m$ 
and $\mu_{0,\infty}=r$, the radius function on $\C^m$, and
$\nu_0=1+r$, so when $i=0$ eqn \eq{qaleeq} becomes
\begin{equation}
\nabla^l\bigl(\psi_0^*(g)-h_0\bigr)=O\bigl(r^{2-2m}(1+r)^{-2-l}\bigr)
\end{equation}
wherever $r>R$ on $\C^m$, for all $l\ge 0$. But this is equivalent to 
the equation \cite[eqn~(3)]{Joyc5} defining ALE K\"ahler metrics. So 
we have proved:

\begin{lem} Suppose $\C^m/G$ has an isolated singularity at\/ $0$. Then 
QALE K\"ahler metrics on $X$ are the same thing as ALE K\"ahler 
metrics on $X$, in the sense of\/~{\rm\cite[Def.~2.3]{Joyc5}}.
\end{lem}

The idea used to prove Proposition \ref{chigiprop} also yields
the following result.

\begin{prop} In the situation of Definition \ref{qaledef}, the 
metrics $g_i$ satisfy $\chi_{\gamma,i}^*(g_{\gamma\cdot i})=g_i$ for 
each\/ $\gamma\in G$ and\/ $i\in I$. Thus $g_i$ is invariant under 
the natural action of\/ $B_i$ on~$Y_i$.
\label{gibiprop}
\end{prop}

Next we show that QALE metrics are made out of other QALE metrics
of lower dimension.

\begin{prop} Let\/ $g_i$ be the K\"ahler metric on the resolution 
$Y_i$ of\/ $W_i/A_i$ in Definition \ref{qaledef}. Then $g_i$ is also 
a QALE metric.
\label{giqaleprop}
\end{prop}

\begin{proof} We shall use the notation defined in the proof of
Proposition \ref{lpprop}. In addition, let $h_j'$ be the Euclidean
metric on $V_j'=V_j\cap W_i$. Let $j\in I'$, so that $i\succeq j$ 
and $V_j=V_i\times V_j'$. Writing the Euclidean metric $h_j$ on $V_j$ 
as $h_i\times h_j'$, eqn \eq{qaleeq} with $i$ replaced by $j$ 
becomes
\begin{equation}
\nabla^l\bigl(\psi_j^*(g)-h_i\times h_j'\times g_j\bigr)=
\sum\begin{Sb}k\in I:j\nsucceq k\end{Sb}
O\bigl(\mu_{j,k}^{d_k}\nu_j^{-2-l}\bigr)
\label{qa1eq}
\end{equation}
on $V_i\times V_j'\times Y_j\,\big\backslash U_j$. Using $\id\times\psi_j'$ 
to pull eqn \eq{qaleeq} from $V_i\times Y_i\,\big\backslash U_i$ back 
to $V_i\times V_j'\times Y_j\,\big\backslash U_j$, and substituting 
$(\id\times\psi_j')^*(\psi_i^*(g))=\psi_j^*(g)$ since
$\psi_i\circ(\id\times\psi_j')=\psi_j$ by \eq{ijjpcdeq}, gives
\begin{equation}
\nabla^l\bigl(\psi_j^*(g)-h_i\times(\psi_j')^*(g_i)\bigr)=
\sum\begin{Sb}k\in I:i\nsucceq k\end{Sb}
O\bigl(\mu_{j,k}^{d_k}\nu_j^{-2-l}\bigr)
\label{qa2eq}
\end{equation}
on $V_i\times V_j'\times Y_j\,\big\backslash U_j$. Hence, subtracting
\eq{qa1eq} and \eq{qa2eq} gives
\begin{equation*}
\nabla^l\bigl(h_i\times(\psi_j')^*(g_j)-h_i\times h_j'\times g_j\bigr)=
\sum\begin{Sb}k\in I:j\nsucceq k\end{Sb}
O\bigl(\mu_{j,k}^{d_k}\nu_j^{-2-l}\bigr).
\end{equation*}
By restricting to $\{v\}\times V_j'\times Y_j$ for $v\in V_i$ and 
taking $v$ to be large, we show that
\begin{equation*}
\nabla^l\bigl((\psi_j')^*(g_i)-h_j'\times g_j\bigr)=
\sum\begin{Sb}k\in I':j\nsucceq k\end{Sb}
O\bigl((\mu'_{j,k})^{d_k}(\nu_j')^{-2-l}\bigr)
\end{equation*}
on $V_j'\times Y_j\,\big\backslash U_j'$, for all $l\ge 0$. 
Thus $g_i$ is a QALE K\"ahler metric on~$Y_i$.
\end{proof}

So local product resolutions and QALE K\"ahler metrics are built 
up by a kind of {\it induction on dimension}. If $X$ is a local 
product resolution with a QALE K\"ahler metric, then the $Y_i$ that 
appear as limits of $X$ are also local product resolutions with 
QALE metrics, but in a lower dimension. This suggests a method 
of proving results about QALE metrics: we assume the result is 
true for all $Y_i$ with $\dim Y_i<\dim X$, and prove it for $X$. 
Then the result holds for all QALE metrics, by induction on 
$\dim X$. We will use this idea several times later on.

Here is the definition of K\"ahler classes and the K\"ahler cone 
on QALE manifolds.

\begin{dfn} Let $X$ be a local product resolution of $\C^m/G$, and $g$ 
a QALE K\"ahler metric on $X$ with K\"ahler form $\omega$. The de Rham 
cohomology class $[\omega]\in H^2(X,\R)$ is called the {\it K\"ahler 
class} of $g$. Define the {\it K\"ahler cone} $\mathcal K$ of $X$ to 
be the set of K\"ahler classes $[\omega]\in H^2(X,\R)$ of QALE K\"ahler 
metrics on~$X$. 
\end{dfn}

By studying the de Rham cohomology of a QALE K\"ahler manifold
one can show that $\mathcal K$ is an open convex cone in 
$H^2(X,\R)$, not containing zero.

\section{Ricci-flat QALE K\"ahler manifolds}

A resolution $(X,\pi)$ of $\C^m/G$ with first Chern class
$c_1(X)=0$ is called a {\it crepant resolution}, as in
Reid \cite{Reid}. A great deal is known about the algebraic
geometry of crepant resolutions, especially in dimensions 2 
and 3. In particular, if $\C^m/G$ has a crepant resolution
then $G\subset{\rm SU}(m)$, and conversely when $m=2,3$ a 
crepant resolution of $\C^m/G$ exists for every finite 
subgroup $G$ of ${\rm SU}(m)$. If $G$ is abelian then all 
the crepant resolutions of $\C^m/G$ can be understood 
explicitly using toric geometry.

Since any Ricci-flat K\"ahler manifold $X$ has $c_1(X)=0$, we 
deduce:

\begin{prop} Suppose $(X,J,g)$ is a Ricci-flat QALE K\"ahler 
manifold asymptotic to $\C^m/G$. Then $G\subset{\rm SU}(m)$ 
and\/ $X$ is a crepant resolution of\/~$\C^m/G$.
\end{prop}

We shall construct Ricci-flat QALE K\"ahler manifolds by
starting with a crepant resolution $X$ of $\C^m/G$, and
constructing a suitable K\"ahler metric on $X$. This is 
why we are interested in crepant resolutions. The proof 
of the following proposition is fairly straightforward 
if you know enough algebraic geometry, and we omit it. 

\begin{prop} Suppose $X$ is a crepant resolution of\/ $\C^m/G$.
Then $X$ is a local product resolution. Moreover, each resolution 
$Y_i$ of\/ $W_i/A_i$ is also a crepant resolution.
\end{prop}

Here is our main result, which generalizes Theorem 3.3 of 
\cite{Joyc5} to the QALE case. Its proof works by first 
proving a version of the Calabi conjecture for QALE manifolds. 

\begin{thm} Let\/ $G$ be a finite subgroup of\/ ${\rm SU}(m)$, and\/ 
$X$ a crepant resolution of\/ $\C^m/G$. Then each K\"ahler class of 
QALE metrics on $X$ contains a unique Ricci-flat QALE K\"ahler metric.
\label{qalerfthm}
\end{thm}

If $X$ is a crepant resolution of $\C^m/G$ then it is a local 
product resolution, by the previous proposition. Also, since 
$G\subset{\rm SU}(m)$ the singular set of $\C^m/G$ has codimension 
$n\ge 2$, which was one of the conditions for $X$ to admit QALE 
metrics in Definition \ref{qaledef}. In practice this means that 
any K\"ahler crepant resolution of $\C^m/G$ also admits QALE 
K\"ahler metrics, and so has a family of Ricci-flat QALE K\"ahler 
metrics by Theorem~\ref{qalerfthm}.

We call a quotient singularity $\C^m/G$ {\it reducible} if 
we can write $\C^m/G=\C^{m_1}/G_1\times\C^{m_2}/G_2$, where 
$G_j\subset{\rm U}(m_j)$ and $m_1,m_2>0$ satisfy $m_1+m_2=m$. 
Otherwise we say $\C^m/G$ is {\it irreducible}. If $\C^m/G$ 
is reducible and $(X,J,g)$ is a Ricci-flat QALE K\"ahler manifold 
asymptotic to $\C^m/G$, then it is easy to see that $(X,J,g)$ 
is a product of lower-dimensional Ricci-flat QALE K\"ahler manifolds.
So we shall restrict our attention to irreducible $\C^m/G$.
In this case we can show:

\begin{thm} Let\/ $(X,J,g)$ be a Ricci-flat QALE K\"ahler manifold 
asymptotic to $\C^m/G$, where $\C^m/G$ is irreducible. If\/ $m\ge 4$ 
is even and\/ $G$ is conjugate to a subgroup of\/ ${\rm Sp}(m/2)$ 
then $\Hol(g)={\rm Sp}(m/2)$, and otherwise~$\Hol(g)={\rm SU}(m)$.
\label{qaleholthm}
\end{thm}

The proofs of Theorems \ref{qalerfthm} and \ref{qaleholthm}  
will be given in \S 7.

\subsection{Examples of QALE manifolds with holonomy ${\rm SU}(m)$}

Here are three examples of QALE manifolds with holonomy ${\rm SU}(m)$,
for $m=3$ and~4. 

\begin{ex} Define $\alpha:\C^3\rightarrow\C^3$ by
\begin{equation*}
\alpha:(z_1,z_2,z_3)\mapsto (-z_1,iz_2,iz_3).
\end{equation*}
Then $\langle\alpha\rangle$ is a subgroup of ${\rm SU}(3)$ isomorphic 
to ${\mathbb Z}_4$. Let $I$ be $\{0,1,\infty\}$, and set
\begin{equation*}
V_0=\C^3,\quad V_1=\bigl\{(z_1,0,0):z_1\in\C\bigr\}
\quad\text{and}\quad V_\infty=\{0\}.
\end{equation*}
Then ${\mathcal L}=\{V_i:i\in I\}$, and the groups $A_i$ are
\begin{equation*}
A_0=\{1\},\quad A_1=\{1,\alpha^2\},\quad\text{and}\quad 
A_\infty=\{1,\alpha,\alpha^2,\alpha^3\}={\mathbb Z}_4.
\end{equation*}

The quotient $\C^3/{\mathbb Z}_4$ has a unique crepant resolution $X$, 
with $b^2(X)=2$. It can be constructed using toric geometry, and admits 
QALE K\"ahler metrics. Thus $X$ has a 2-parameter family of Ricci-flat 
QALE K\"ahler metrics $g$ by Theorem \ref{qalerfthm}, which have 
holonomy ${\rm SU}(3)$ by Theorem~\ref{qaleholthm}.

We can describe the asymptotic behaviour of these Ricci-flat 
metrics $g$ quite simply. The vector space $W_1$ is
$\bigl\{(0,z_2,z_3):z_2,z_3\in\C\bigr\}$, and $\alpha^2\in A_1$ 
acts as $-1$ on $W_1$. Thus $W_1/A_1\cong\C^2/\{\pm1\}$. Let $Y_1$ 
be the blow-up of $\C^2/\{\pm1\}$ at 0. Then $Y_1$ carries a 
1-parameter family of ALE metrics with holonomy SU(2), given 
explicitly by Eguchi and Hanson \cite{EgHa}. Let $g_1$ be an 
Eguchi-Hanson metric on $Y_1$, and $h_1$ the Euclidean metric 
on $V_1\cong\C$. Then $h_1\times g_1$ is a Ricci-flat metric 
on~$V_1\times Y_1$. 

Now $\alpha$ acts on $V_1\times Y_1$ with $\alpha^2=1$, and 
the fixed points of $\alpha$ are a copy of $\mathbb{CP}^1$ in 
$V_1\times Y_1$. Thus $(V_1\times Y_1)/\langle\alpha\rangle$ is 
a complex orbifold, and $X$ is its unique crepant resolution. The 
metric $h_1\times g_1$ is preserved by $\alpha$ and descends to 
$(V_1\times Y_1)/\langle\alpha\rangle$. Identifying $X$ with 
$(V_1\times Y_1)/\langle\alpha\rangle$ outside $\pi^{-1}(0)$, the 
asymptotic conditions \eq{qaleeq} on our QALE metric $g$ become
\begin{equation*}
\nabla^l(g-h_1\times g_1)=O\bigl(\pi^*(r)^{-4}(1+\pi^*(s))^{-2-l}\bigr),
\end{equation*}
where $r,s:\C^3/{\mathbb Z}_4\rightarrow[0,\infty)$ are the 
distances in $\C^3/{\mathbb Z}_4$ to 0 and to the singular set 
$S=\bigl\{\pm(z_1,0,0):z_1\in\C\bigr\}$ respectively. Thus at 
infinity $g$ is asymptotic to $h_1\times g_1$, a metric we can 
write down explicitly in coordinates.
\label{c3/z4ex}
\end{ex}

\begin{ex} Define $\alpha,\beta:\C^3\rightarrow\C^3$ by
\begin{equation*}
\alpha:(z_1,z_2,z_3)\mapsto(z_1,-z_2,-z_3),\quad
\beta:(z_1,z_2,z_3)\mapsto(-z_1,z_2,-z_3).
\end{equation*}
Then $\langle\alpha,\beta\rangle$ is a subgroup of ${\rm SU}(3)$ 
isomorphic to ${\mathbb Z}_2^2$. Let $I$ be $\{0,1,2,3,\infty\}$, 
and set
\begin{equation*}
\begin{alignedat}{3} 
V_0&=\C^3,\qquad V_\infty=\{0\},&\quad 
V_1&=\bigl\{(z_1,0,0):z_1\in\C\bigr\},\\
V_2&=\bigl\{(0,z_2,0):z_2\in\C\bigr\}\quad\text{and}&\quad
V_3&=\bigl\{(0,0,z_3):z_3\in\C\bigr\}.
\end{alignedat}
\end{equation*}
Then ${\mathcal L}=\{V_i:i\in I\}$, and the groups $A_i$ are
\begin{equation*}
A_0=\{1\},\quad 
A_1=\{1,\alpha\},\quad 
A_2=\{1,\beta\},\quad
A_3=\{1,\alpha\beta\},\quad
A_\infty={\mathbb Z}_2^2.
\end{equation*}

The quotient $\C^3/{\mathbb Z}_2^2$ has four distinct crepant 
resolutions $X_1,\ldots,X_4$, with $b^2(X_j)=3$. They can be 
constructed explicitly using toric geometry, and all admit QALE 
K\"ahler metrics. Thus $X_1,\ldots,X_4$ carry 3-parameter families 
of Ricci-flat QALE K\"ahler metrics by Theorem \ref{qalerfthm}, 
which have holonomy ${\rm SU}(3)$ by Theorem~\ref{qaleholthm}.
\label{c3/z22ex}
\end{ex}

\begin{ex} Define $\alpha,\beta,\gamma:\C^4\rightarrow\C^4$ by
\begin{equation*}
\begin{split}
\alpha:&(z_1,\ldots,z_4)\mapsto (-z_1,-z_2,z_3,z_4),\\
\beta:&(z_1,\ldots,z_4)\mapsto (z_1,-z_2,-z_3,z_4)\\
\text{and}\qquad
\gamma:&(z_1,\ldots,z_4)\mapsto (z_1,z_2,-z_3,-z_4).
\end{split}
\end{equation*}
Then $\langle\alpha,\beta,\gamma\rangle\cong{\mathbb Z}_2^3$ 
is a subgroup of ${\rm SU}(4)$. Let $I$ be $\{0,\infty\}\cup
\{jk:j,k=1,\ldots,4$, $j<k\}$, an 8-element set. Then 
${\mathcal L}=\{V_i:i\in I\}$, where
\begin{equation*}
V_0=\C^4,\quad V_\infty=\{1\},\quad
V_{jk}=\bigl\{(z_1,\ldots,z_4)\in\C^4:z_j=z_k=0\bigr\}.
\end{equation*}
The quotient $\C^4/{\mathbb Z}_2^3$ has 48 distinct crepant resolutions 
$X_1,\ldots,X_{48}$, with $b^2(X_j)=6$. They can be constructed 
explicitly using toric geometry, and all admit QALE K\"ahler metrics. 
Thus $X_1,\ldots,X_{48}$ carry 6-parameter families of Ricci-flat 
QALE K\"ahler metrics by Theorem \ref{qalerfthm}, which have holonomy 
${\rm SU}(4)$ by Theorem~\ref{qaleholthm}.
\label{c4/z23ex}
\end{ex}

In each of these examples we know by Theorem \ref{qalerfthm} that 
QALE metrics with holonomy ${\rm SU}(m)$ exist on the resolutions $X$ 
of $\C^m/G$, but we are unable to write these metrics down explicitly
in coordinates. The author believes that QALE metrics with holonomy 
${\rm SU}(m)$ for $m\ge 3$ are not algebraic objects, and it is not 
possible to give an explicit formula for these metrics in coordinates.

\subsection{Examples of QALE manifolds with holonomy ${\rm Sp}(m)$}

It can be shown \cite[\S 9.3]{Joyc4} that if $X$ is a QALE
manifold with holonomy ${\rm Sp}(m)$ for $m\ge 2$ asymptotic
to $\C^{2m}/G$, then $G$ is nonabelian. So to find examples
we must consider nonabelian groups.

\begin{ex} Define $\alpha,\beta:\C^4\rightarrow\C^4$ by
\begin{equation*}
\begin{aligned}
\alpha:&(z_1,\ldots,z_4)\mapsto (e^{2\pi i/3}z_1,e^{4\pi i/3}z_2,
e^{4\pi i/3}z_3,e^{2\pi i/3}z_4),\\
\beta:&(z_1,\ldots,z_4)\mapsto (z_3,z_4,z_1,z_2).
\end{aligned}
\end{equation*}
Then $G=\langle\alpha,\beta\rangle$ is a nonabelian subgroup of Sp(2) 
of order 6 isomorphic to the symmetric group $S_3$, that preserves 
the complex symplectic form~$\d z_1\wedge\d z_2+\d z_3\wedge\d z_4$.

Let $S$ be the singular set of $\C^4/G$, and define a map
$\phi:\C^4/G\,\big\backslash S\rightarrow\mathbb{CP}^4$ by
\begin{equation}
\begin{split}
\phi\bigl((z_1,\ldots,z_4)G\bigr)=[&z_1z_2-z_3z_4,z_1^3-z_3^3,\\
&z_1^2z_4-z_2z_3^2,
z_1z_4^2-z_2^2z_3,z_2^3-z_4^3].
\end{split}
\end{equation}
The five polynomials in $z_1,\ldots,z_4$ given here are invariant 
under $\alpha$ and change sign under $\beta$, and they are all zero 
if and only if $(z_1,z_2,z_3,z_4)G$ lies in $S$. Let $X$ be the 
closure of the graph of $\phi$ in $\C^4/G\times\mathbb{CP}^4$, and 
$\pi:X\rightarrow\C^4/G$ the natural projection. Then a careful 
analysis shows that $X$ is actually nonsingular, and $(X,\pi)$ is a 
{\it crepant resolution} of~$\C^4/G$.

Let $[x_0,x_1,x_2]$ be homogeneous coordinates on the {\it 
weighted projective space} $\mathbb{CP}^2_{3,1,1}$. Define a 
map $\psi:\mathbb{CP}^2_{3,1,1}\rightarrow\mathbb{CP}^4$ by
\begin{equation}
\psi\bigl([x_0,x_1,x_2]\bigr)=[x_0,x_1^3,x_1^2x_2,x_1x_2^2,x_2^3].
\end{equation}
Then $\psi$ is injective, and it can be shown that $\pi^{-1}(0)$ is 
the subset $\{0\}\times\mathop{\rm Im}\psi$ in $\C^4/G\times\mathbb{CP}^4$. 
Thus $\pi^{-1}(0)$ is isomorphic to the (singular) weighted projective 
space $\mathbb{CP}^2_{3,1,1}$. Since $X$ retracts onto $\pi^{-1}(0)$, it 
follows that $b^2(X)=1$, $b^4(X)=1$ and $b^6(X)=0$. Clearly $X$ is a 
quasi-projective variety, and therefore K\"ahler. Theorem \ref{qalerfthm} 
shows that $X$ has a 1-parameter family of Ricci-flat QALE K\"ahler 
metrics, which have holonomy Sp(2) by Theorem~\ref{qaleholthm}.
\label{c4/s3ex}
\end{ex}

In fact Example \ref{c4/s3ex} generalizes to give an action 
of the symmetric group $S_{m+1}$ on $\C^{2m}$, and a crepant 
resolution $X_m$ of $\C^{2m}/S_{m+1}$ carrying QALE K\"ahler 
metrics with holonomy ${\rm Sp}(m)$. To prove this we adapt 
an idea of Beauville, and regard the {\it Hilbert scheme} or 
{\it Douady space} ${\mathbin{\rm Hilb}}^{m+1}\C^2$ of $m+1$ 
points in $\C^2$ as a crepant resolution of $\C^{2m+2}/S_{m+1}$. 
Further details are given in~\cite[\S 6-7]{Beau}.

\begin{ex} Let $\C^4$ have coordinates $(z_1,z_2,z_3,z_4)$. Let $H$ 
be a subgroup of SU(2). Then $H\times H$ acts on $\C^4$, with the 
first $H$ acting only on the coordinates $(z_1,z_2)$ and the second 
$H$ acting only on $(z_3,z_4)$. Define $\alpha:\C^4\rightarrow\C^4$ by
$\alpha:(z_1,\ldots,z_4)\mapsto (z_3,z_4,z_1,z_2)$. Let $G$ be the 
subgroup of Sp(2) generated by $H\times H$ and $\alpha$. Then $G$ is a 
semidirect product ${\mathbb Z}_2\ltimes(H\times H)$. We can define a 
crepant resolution $X$ of $\C^4/G$ as follows.

Let $Y$ be the unique crepant resolution of $\C^2/H$. Then $Y\times Y$ 
is a crepant resolution of $\C^4/(H\times H)$, and the action of $\alpha$ 
lifts to $Y\times Y$ by $\alpha:(y_1,y_2)\mapsto(y_2,y_1)$ in the obvious 
way. The singular set of the quotient $(Y\times Y)/\langle\alpha\rangle$ 
is the `diagonal' $\Delta_Y=\bigl\{(y,y):y\in Y\bigr\}$, and each singular 
point is modelled locally on $\C^2\times(\C^2/\{\pm1\})$. Let $X$ be the 
blow-up of $(Y\times Y)/\langle\alpha\rangle$ along the diagonal $\Delta_Y$. 
Then $X$ is nonsingular, and is a crepant resolution of~$\C^4/G$. 

Since we understand $Y$ very well, it is easy to compute the
Betti numbers of $X$. For example, in the case $H={\mathbb Z}_k$ we 
have $b^2(X)=k$, $b^4(X)={1\over 2}(k+2)(k-1)$ and $b^6(X)=0$. Again,
each such resolution $X$ has a family of Ricci-flat QALE K\"ahler 
metrics by Theorem \ref{qalerfthm}, which have holonomy Sp(2) 
by Theorem~\ref{qaleholthm}.
\end{ex}

We claimed in \S 3.1 that QALE metrics with holonomy
${\rm SU}(m)$ for $m\ge 3$ are nonalgebraic objects, and cannot be
explicitly written down in coordinates. However, for QALE metrics 
with holonomy ${\rm Sp}(m)$, the reverse is true. The author has a 
proof that {\it every} QALE metric with holonomy ${\rm Sp}(m)$ has 
an algebraic description. This proof uses the theory of
{\it hypercomplex algebraic geometry}, which is described in
Joyce~\cite{Joyc3}. 

It seems likely that QALE manifolds with holonomy ${\rm Sp}(m)$ can 
be explicitly constructed using the hyperk\"ahler quotient, as in 
Kronheimer's construction of ALE manifolds with holonomy Sp(1),
\cite{Kron}. However, at present the author has no proof of this, 
nor any explicit examples.

\section{K\"ahler potentials on QALE K\"ahler manifolds}

If $g,g'$ are two QALE K\"ahler metrics in the same K\"ahler class 
on $X$, with K\"ahler forms $\omega,\omega'$, then we expect that 
$\omega'=\omega+\d\d^c\phi$ for some function $\phi$ on $X$. (Here we 
use the notation that $\d^c\phi=i(\overline\partial-\partial)\phi$, so 
that $\d^c$ is a real operator with $\d\d^c=2i\partial\overline\partial$.) 
Conversely, if $g$ is a QALE K\"ahler metric on $X$ with K\"ahler form 
$\omega$, then we can try to define other QALE metrics $g'$ on $X$ with 
K\"ahler forms $\omega'=\omega+\d\d^c\phi$ for suitable functions $\phi$ 
on $X$. In this section we will study the properties of such functions 
$\phi$ in detail.

\begin{dfn} Let $(X,J,g)$ be a QALE K\"ahler manifold asymptotic to
$\C^m/G$, and use the notation of Definition \ref{qaledef}. We 
say that a smooth real function $\phi$ on $X$ is of {\it K\"ahler 
potential type} if for each $i\in I$ there exists a smooth 
real function $\phi_i$ on $Y_i$, with $\phi_0=0$, such that 
\begin{equation}
\nabla^l\bigl(\psi_i^*(\phi)-\phi_i\bigr)=
\sum\begin{Sb}j\in I:i\nsucceq j\end{Sb}
O\bigl(\mu_{i,j}^{d_j}\nu_i^{-l}\bigr)
\label{qphieq}
\end{equation}
on $V_i\times Y_i\,\big\backslash U_i$, for all $l\ge 0$. Here we 
identify $\phi_i$ with its pull-back to~$V_i\times Y_i$.
\label{qktypedef}
\end{dfn}

By comparing \eq{qphieq} with \eq{qaleeq}, we immediately deduce:

\begin{prop} Let\/ $(X,J,g)$ be a QALE K\"ahler manifold, $\omega$ 
the K\"ahler form of $g$, and\/ $\phi$ a function of K\"ahler 
potential type on $X$. Suppose $\omega'=\omega+\d\d^c\phi$ is a 
positive $(1,1)$-form on $X$. Then the K\"ahler metric $g'$ on $X$ 
with K\"ahler form $\omega'$ is~QALE.
\label{kpotqaleprop}
\end{prop}

Here is a converse to this proposition. It can be proved by 
combining the method of \cite[Th.~5.4]{Joyc5} with the results 
on the Laplacian on QALE manifolds that we will prove in \S 5, 
extended along the lines sketched in \cite[\S 9.7]{Joyc4}. We
will not give a proof here.

\begin{thm} Let\/ $X$ be a local product resolution of\/ 
$\C^m/G$, and suppose that\/ $g,g'$ are QALE metrics on $X$ in 
the same K\"ahler class, with K\"ahler forms $\omega,\omega'$. 
Then there is a unique function $\phi$ of K\"ahler potential 
type on $X$ such that\/~$\omega'=\omega+\d\d^c\phi$.
\label{kpotqalethm}
\end{thm}

Proposition \ref{kpotqaleprop} and Theorem \ref{kpotqalethm} show
that functions of K\"ahler potential type are the natural class of
functions to use as K\"ahler potentials on QALE K\"ahler manifolds. 
The next result can be proved by following the proofs of Propositions 
\ref{gibiprop} and~\ref{giqaleprop}.

\begin{prop} Let\/ $(X,J,g)$ be a QALE K\"ahler manifold, and\/ 
$\phi$ a function of K\"ahler potential type on $X$. Then the 
functions $\phi_i$ on $Y_i$ introduced in Definition \ref{qktypedef}
satisfy the following conditions:
\begin{itemize}
\item[{\rm(i)}] For each\/ $\gamma\in G$ and\/ $i\in I$ we have
$\chi_{\gamma,i}^*(\phi_{\gamma\cdot i})=\phi_i$. Hence
$\phi_i$ is invariant under the action of\/ $B_i$ on~$Y_i$.
\item[{\rm(ii)}] Whenever $i,j\in I$ with $i\succeq j$, then on 
$V_j'\times Y_j'\,\big\backslash U_j'$ we have 
\begin{equation*}
\nabla^l\bigl((\psi_j')^*(\phi_i)-\phi_j\bigr)=
\sum\begin{Sb}k\in I:\\
i\succeq k,\;j\nsucceq k\end{Sb}
O\bigl((\mu_{j,k}')^{d_k}(\nu_j')^{-l}\bigr)
\end{equation*}
for all\/ $l\ge 0$, using the notation of Propositions 
\ref{lpprop} and \ref{giqaleprop}. Thus the functions 
$\phi_i$ on $Y_i$ are also of K\"ahler potential type.
\end{itemize}
\label{qkpotcondprop}
\end{prop}

The proposition shows that (i) and (ii) are {\it necessary}
conditions for a set of functions $\phi_i$ to be associated 
to a function $\phi$ of K\"ahler potential type. In fact they are 
also {\it sufficient}.

\begin{thm} Let\/ $(X,J,g)$ be a QALE K\"ahler manifold, and suppose 
that for each\/ $i\in I\setminus\{\infty\}$ there is a smooth function 
$\phi_i$ on $Y_i$, with\/ $\phi_0=0$, such that the $\phi_i$ satisfy 
conditions {\rm(i)} and\/ {\rm(ii)} of Proposition \ref{qkpotcondprop}. 
Then there exists a function $\phi$ of K\"ahler potential type on $X$
asymptotic to these $\phi_i$ for all\/~$i\in I\setminus\{\infty\}$.
\label{phiextthm}
\end{thm}

\begin{proof} Let $\eta:[0,\infty)\rightarrow[0,1]$ be smooth 
with $\eta(x)=0$ for $x\le R$ and $\eta(x)=1$ for $x\ge 2R$,
where $R>0$ is the constant in Definition \ref{locproddef}.
For each $i\in I\setminus\{\infty\}$, define a smooth function 
$\Phi_i$ on $V_i\times Y_i$ by
\begin{equation}
\Phi_i(v,y)=\phi_i(y)\cdot 
\prod\begin{Sb}j\in I:i\nsucceq j\end{Sb}\eta\bigl(\mu_{i,j}(v,y)\bigr).
\label{phiipdefeq}
\end{equation}
The idea here is that $\Phi_i=\phi_i$ at distance at least $2R$ from 
the pull-back of $S_i$ in $V_i\times Y_i$, that $\Phi_i=0$ at distance 
no more than $R$ from the pull-back of $S_i$, and that between 
distances $R$ and $2R$ we join the two possibilities $\Phi_i=\phi_i$ 
and $\Phi_i=0$ smoothly together using a partition of unity. 
Note that $\Phi_i\equiv 0$ in~$U_i$.

For $i\in I\setminus\{\infty\}$, let $k_i$ be integers satisfying
\begin{equation}
\sum\begin{Sb}i\in I\setminus\{\infty\}:i\succeq j\end{Sb}
k_i=1\qquad\text{for each $j\in I\setminus\{\infty\}$.}
\label{qalekidef}
\end{equation}
It can be shown that these equations have a unique solution 
$\{k_i\}$. Now define
\begin{equation}
\phi(x)=\sum_{i\in I\setminus\{\infty\}}{k_i\md{A_i}\over\md{G}}
\sum\begin{Sb}(v,y)\in V_i\times Y_i\setminus U_i:\\
\psi_i(v,y)=x\end{Sb}\Phi_i(v,y).
\label{qalephiieq}
\end{equation}
As $\Phi_i\equiv 0$ in $U_i$, we see that $\phi$ is smooth. It 
turns out that this $\phi$ is of K\"ahler potential type on $X$
and asymptotic to $\phi_i$ for all $i\in I\setminus\{\infty\}$, 
so that it satisfies the conditions of the theorem. The proof 
of this is rather complicated, and we will not give it in full. 
Instead, we will explain the important points under the 
simplifying assumption that $N(V_i)=G$ for all~$i\in I$.

In this case $B_i=G/A_i$, so $\md{B_i}=\md{G}/\md{A_i}$. Thus 
we may rewrite \eq{qalephiieq} as
\begin{equation}
\phi(x)=\sum_{i\in I\setminus\{\infty\}}k_i\Phi_i'(x),
\quad\text{where}\quad
\Phi_i'(x)={1\over\md{B_i}}
\sum\begin{Sb}(v,y)\in V_i\times Y_i\setminus U_i:\\ 
\psi_i(v,y)=x\end{Sb}\Phi_i(v,y).
\label{qalephijeq}
\end{equation}
Now $B_i$ acts on $Y_i$ and $V_i\times Y_i$, and $\phi_i$ is 
$B_i$-invariant by condition (i). Thus $\Phi_i$ is also 
$B_i$-invariant, and pushes down to a function on $(V_i\times Y_i)/B_i$. 
From \S 1, $\psi_i$ induces an isomorphism between 
$(V_i\times Y_i\setminus U_i)/B_i$ and $X\setminus\pi^{-1}(T_i/B_i)$. 
In fact $\Phi_i'$ is the push-forward of $\Phi_i$ under this isomorphism. 
The factor $1/\md{B_i}$ in \eq{qalephijeq} ensures this, because 
each point of $X\setminus\pi^{-1}(T_i/B_i)$ pulls back to an orbit of 
$B_i$ in $V_i\times Y_i$, consisting of $\md{B_i}$ points. 

Thus, abusing notation a little, we can write $\phi=\sum_{i\ne\infty}
k_i\Phi_i$, because the factor $\md{A_i}/\md{G}$ in \eq{qalephiieq} 
compensates for the fact that each generic point in $X$ pulls back 
to $\md{G}/\md{A_i}$ points in $V_i\times Y_i$. But $\Phi_i=\phi_i$
away from $U_i$, and so away from $\pi^{-1}(S)$ in $X$ we have
$\phi=\sum_{i\ne\infty}k_i\phi_i$, again by an abuse of notation.

We must prove that $\phi$ satisfies \eq{qphieq}. One way to 
interpret this is to say that $X$ is divided roughly into 
overlapping regions corresponding to $j\in I$, where in the 
$j^{\rm th}$ region $X$ is locally isomorphic to $V_j\times Y_j$ 
and $\phi\approx\phi_j$. Now on the $j^{\rm th}$ region we 
have $\phi_i\approx\phi_j$ if $i\succeq j$ and $\phi_i\approx 0$ 
if $i\nsucceq j$, by condition (ii) of Proposition 
\ref{kpotqaleprop}. Therefore on the $j^{\rm th}$ region we have
\begin{equation*}
\phi\approx\sum\begin{Sb}i\in I\setminus\{\infty\}:i\succeq j\end{Sb}
k_i\phi_j=\phi_j,
\end{equation*}
by \eq{qalekidef}, as we want. This idea can be used to show 
that \eq{qphieq} holds away from $\pi^{-1}(S)$ in~$X$.

It remains to consider the parts of $X$ near $\pi^{-1}(S)$.
Within distance $2R$ of $\pi^{-1}(S)$ in $X$, the functions
$\eta(\mu_{i,j})$ in \eq{phiipdefeq} are not all equal to 1, 
and so we do not have $\Phi_i=\phi_i$ for all $i$. The 
reason for introducing the $\eta(\mu_{i,j})$ is that the 
push-forward of $\phi_i$ to $X\setminus\pi^{-1}(S_i/B_i)$ may
not extend smoothly to $X$. So we modify $\phi_i$ to get 
$\Phi_i$, whose push-forward is zero near $\pi^{-1}(S_i/B_i)$ 
and does extend smoothly to $X$. But we then have to make sure 
that the `errors' due to the $\eta(\mu_{i,j})$ are within the 
bounds allowed by \eq{qphieq}. One can show using \eq{qalekidef} 
that this is always so, and the proof is complete.
\end{proof}

The reason for excluding $i=\infty$ in this theorem is that in 
general $X=V_\infty\times Y_\infty$ and $\phi=\phi_\infty$, and 
so if we assumed that $\phi_\infty$ existed we could just take 
$\phi=\phi_\infty$ and there would be nothing to prove. Although 
the proof works by writing down $\phi$ explicitly, there are in 
fact many suitable $\phi$, since if we add to $\phi$ any smooth 
function on $Y_\infty$ with sufficiently fast decay at infinity, 
the resulting function also satisfies the theorem. The proof above 
will play an important part in the construction of Ricci-flat QALE 
metrics in \S 7, which is why we covered it in some detail.

\begin{prop} In the situation of Theorem \ref{phiextthm}, let\/
$\omega$ be the K\"ahler form of\/ $g$ and\/ $\omega_i$ the 
K\"ahler form of\/ $g_i$, and suppose $\omega_i+\d\d^c\phi_i$ is a 
positive $(1,1)$-form on $Y_i$ for each $i\in I\setminus\{\infty\}$. 
Then we can choose the function $\phi$ such that\/ $\omega'=
\omega+\d\d^c\phi$ is a positive $(1,1)$-form on $X$. Thus the 
K\"ahler metric $g'$ on $X$ with K\"ahler form $\omega'$ is QALE, 
by Proposition~\ref{kpotqaleprop}.
\label{phiextprop}
\end{prop}

\begin{proof} Let $\phi$ be as in the proof of Theorem 
\ref{phiextthm}. Then $\omega+\d\d^c\phi$ is positive outside a 
compact set in $X$, because at large distances from $\pi^{-1}(0)$ 
we have $\phi\approx\phi_i$ for some $i\ne\infty$, and so 
$\omega+\d\d^c\phi$ is positive because $\omega_i+\d\d^c\phi_i$ is 
positive. Let $\hat R>0$, and choose a smooth function 
$\eta:X\rightarrow [0,1]$ such that $\eta=0$ at distance less than 
$\hat R$ from $\pi^{-1}(0)$ in $X$, and $\eta=1$ at distance more than 
$2\hat R$ from $\pi^{-1}(0)$ in $X$, and $\nabla\eta=O(\hat R^{-1})$, 
$\nabla^2\eta=O(\hat R^{-2})$. For large $\hat R$ it turns out that 
$\eta\phi$ satisfies the conditions of the proposition. 
\end{proof}

One moral of this proposition is that QALE K\"ahler metrics on $X$ 
are actually very abundant. That is, the asymptotic conditions
on QALE metrics are not so restrictive that they admit few
solutions. Rather, given any set of K\"ahler metrics $g_i$ on $Y_i$
for $i\in I\setminus\{\infty\}$ satisfying the obvious consistency 
conditions, we expect to find many QALE K\"ahler metrics $g$ on 
$X$ asymptotic to these~$g_i$.

Finally, we give an analogue of \cite[Prop.~5.5]{Joyc5} for 
QALE manifolds.

\begin{thm} Let\/ $X$ be a local product resolution of\/ $\C^m/G$ that 
admits QALE K\"ahler metrics. Then in each K\"ahler class there exists
a QALE K\"ahler metric $g'$ on $X$ such that for each\/ $i\in I$ we have
$\psi_i^*(g')=h_i\times g_i'$ on the set\/ $V_i\times Y_i\setminus U_i$, 
where $g_i'$ is a K\"ahler metric on $Y_i$, and\/ $U_i$ is defined as 
in Definition \ref{locproddef} using a constant\/ $R>0$, which may 
depend on the K\"ahler class.
\label{qalekthm}
\end{thm}

Taking $i=0$ gives $g'=\pi^*(h_0)$ on 
$\bigl\{x\in X:d(\pi(x),S)>R\bigr\}$, where $h_0$ is the 
Euclidean metric, $S$ the singular set and $d(\,,\,)$ the 
distance on $\C^m/G$. Thus $g'$ is flat outside a fixed 
distance from the exceptional set of the resolution $X$. 
Also, $g'$ is a product metric $h_i\times g_i'$ on the parts 
of $X$ modelled on $V_i\times Y_i$, except within a fixed 
distance of some other component of the exceptional set.

Here is a sketch of the proof of this theorem; we leave the details 
as an exercise for the reader. The basic idea is to start with a 
QALE K\"ahler manifold $(X,J,g)$ with K\"ahler form $\omega$, and 
to find a function $\phi$ of K\"ahler potential type on $X$ such 
that $\omega'=\omega+\d\d^c\phi$ is a positive (1,1)-form, and the 
QALE K\"ahler metric $g'$ with K\"ahler form $\omega'$ satisfies 
the conditions of the theorem. We construct $\phi$ by induction 
on~$m=\dim X$.

The inductive step of the proof works as follows. Assume by
induction that such a $\phi$ exists when $\dim X<m$, and
let $(X,J,g)$ be QALE K\"ahler manifold asymptotic to $\C^m/G$.
The inductive hypothesis shows that for each $i\ne\infty$ in 
$I$ there is a function $\phi_i$ on $Y_i$ such that 
$\omega_i'=\omega_i+\d\d^c\phi_i$ is positive on $Y_i$, and the 
corresponding K\"ahler metric $g_i'$ on $Y_i$ satisfies the theorem. 

We can also arrange that these functions $\phi_i$ for 
$i\ne\infty$ satisfy conditions (i) and (ii) of Proposition
\ref{qkpotcondprop}. Then using the ideas of Theorem 
\ref{phiextthm} and Proposition \ref{phiextprop}, we
construct a function $\phi$ on $X$ with the properties
we need, that is asymptotic to $\phi_i$ on $Y_i$ for
each $i\ne\infty$. To prove such a $\phi$ exists requires
the analytical results of \S 5 and also some
discussion of the de Rham cohomology $H^2(X,\R)$, along
the lines of the proof of~\cite[Prop.~5.5]{Joyc5}.

\section{Analysis on QALE K\"ahler manifolds}

In our previous paper \cite[\S 4]{Joyc5} we described the
{\it weighted H\"older spaces} $C^{k,\alpha}_\beta(X)$
over an ALE manifold $X$. These are Banach spaces of functions 
on $X$ defined using powers of a {\it radius function} $\rho$ 
on $X$. They are important tools in solving analysis problems 
on ALE manifolds, because elliptic operators such as the 
Laplacian $\Delta$ have good regularity properties on 
these spaces. 

In this section we develop a similar theory of analysis on QALE 
manifolds. We shall define weighted H\"older spaces of functions 
on QALE K\"ahler manifolds $X$, and study the action of the 
Laplacian on them. We begin by defining functions $\rho,\sigma$ 
on $X$ which are analogues of the radius functions used 
in~\cite{Joyc5}.

\begin{dfn} Let $G$ be a finite subgroup of ${\rm U}(m)$, and $S$ 
the singular set of $\C^m/G$. Define continuous functions 
$r:\C^m/G\rightarrow[0,\infty)$ and $s:\C^m/G\rightarrow[0,\infty)$ by 
$r(x)=d(x,0)$ and $s(x)=d(x,S)$, where $d(\,,\,)$ is the distance on 
$\C^m/G$. Let $(X,\pi)$ be a local product resolution of $\C^m/G$, and
$g$ a QALE K\"ahler metric on $X$. We say that $(\rho,\sigma)$ is a {\it 
pair of radius functions on} $X$ if $\rho,\sigma:X\rightarrow[1,\infty)$ 
are smooth functions such that $\rho\ge\sigma$ and for some $K>0$ we have
\begin{alignat}{7}
\pi^*(r)+1&\le\rho\le&\,\pi^*(r)+2,&\quad
\md{\nabla\rho}&\le K,\quad&\text{and}&\quad
\md{\nabla^2\rho}&\le K\rho^{-1},
\label{rhoesteq}\\
{\textstyle{1\over 2}}\pi^*(s)+1&\le\sigma\le&\,\pi^*(s)+2,&\quad
\md{\nabla\sigma}&\le K,\quad&\text{and}&\quad
\md{\nabla^2\sigma}&\le K\sigma^{-1}.
\label{siesteq}
\end{alignat}
A pair of radius functions $(\rho,\sigma)$ exists for every 
QALE K\"ahler manifold~$(X,J,g)$.
\label{qalerfdef}
\end{dfn}

Here $\rho$ and $\sigma$ are smoothed versions of $\pi^*(r)$ 
and $\pi^*(s)$, adjusted to ensure that $\rho\ge 1$ and
$\sigma\ge 1$. The reason why we have $\pi^*(r)+1\le\rho$
but ${1\over 2}\pi^*(s)+1\le\sigma$ is that $s$ is a rather less
smooth function on $\C^m/G$ than $r$, and so to make a
smoothed version $\sigma$ with $\nabla^2\sigma$ small we must
allow a greater difference between $\sigma$ and $\pi^*(s)$
than between $\rho$ and~$\pi^*(r)$.

The author believes that the following is a natural 
definition of weighted H\"older spaces on QALE manifolds. 

\begin{dfn} Let $(X,J,g)$ be a QALE K\"ahler manifold and $(\rho,\sigma)$ 
a pair of radius functions on $X$. For $\beta,\gamma\in\R$ and $k$ a 
nonnegative integer, define $C^k_{\beta,\gamma}(X)$ to be the 
space of continuous functions $f$ on $X$ with $k$ continuous 
derivatives, such that $\rho^{-\beta}\sigma^{j-\gamma}\bmd{\nabla^jf}$ 
is bounded on $X$ for $j=0,\ldots,k$. Define the norm 
$\nm{\,.\,}_{\smash{C^k_{\beta,\gamma}}}$ on $C^k_{\beta,\gamma}(X)$ by 
\begin{equation*}
\nm{f}_{\smash{C^k_{\beta,\gamma}}}=\sum_{j=0}^k\sup_X
\bmd{\rho^{-\beta}\sigma^{j-\gamma}\nabla^jf}.
\end{equation*}
Let $\delta(g)$ be the injectivity radius of $g$, and write $d(x,y)$ 
for the distance between $x,y$ in $X$. For $T$ a tensor field on $X$ 
and $\alpha,\beta,\gamma\in\R$, define
\begin{equation}
\begin{split}
\bigl[T\bigr]_{\alpha,\beta,\gamma}=
\sup\begin{Sb}x\ne y\in X\\ d(x,y)<\delta(g)\end{Sb}
\biggl[
&\min\bigl(\rho(x),\rho(y)\bigr)^{-\beta}\cdot\\[-10pt]
&\min\bigl(\sigma(x),\sigma(y)\bigr)^{-\gamma}\cdot
{\bmd{T(x)-T(y)}\over d(x,y)^\alpha}\,\biggr].
\end{split}
\label{qaleholdereq}
\end{equation}
Here we interpret $\md{T(x)-T(y)}$ using parallel translation along 
the unique geodesic of length $d(x,y)$ joining $x$ and $y$. For 
$\beta,\gamma\in\R$, $k$ a nonnegative integer, and $\alpha\in(0,1)$, 
define the {\it weighted H\"older space} $C^{k,\alpha}_{\beta,\gamma}(X)$ 
to be the set of $f\in C^k_{\beta,\gamma}(X)$ for which the norm
\begin{equation}
\bnm{f}_{C^{k,\alpha}_{\beta,\gamma}}=\bnm{f}_{C^k_{\beta,\gamma}}
+\bigl[\nabla^k f\bigr]_{\alpha,\beta,\gamma-k-\alpha}
\label{qaleholdnormeq}
\end{equation}
is finite. Define $C^\infty_{\beta,\gamma}(X)$ to be the 
intersection of the $C^k_{\beta,\gamma}(X)$ for all $k\ge 0$. Both 
$C^k_{\beta,\gamma}(X)$ and $C^{k,\alpha}_{\beta,\gamma}(X)$ are 
Banach spaces, but $C^\infty_{\beta,\gamma}(X)$ is not a Banach space.
\label{qholderdef}
\end{dfn}

This generalizes the definition of weighted H\"older spaces
$C^{k,\alpha}_\beta(X)$ on ALE manifolds, \cite[Def.~4.2]{Joyc5}.
The reasoning leading up to the definition is discussed in 
\cite[\S 9.5]{Joyc4}, where we also define weighted Sobolev 
spaces on QALE manifolds. The basic idea is that if 
$C^k_{\beta,\gamma}(X)$ then 
$\nabla^jf=O(\rho^\beta\sigma^{\gamma-j})$ for $j\le k$. If 
$\C^m/G$ has an isolated singularity, so that $X$ is an ALE 
manifold, then $\rho=\sigma$ and $C^{k,\alpha}_{\beta,\gamma}(X)$ 
agrees with $C^{k,\alpha}_{\beta+\gamma}(X)$, the usual weighted 
H\"older space on an ALE manifold.

In \cite[Th.~4.5(a)]{Joyc5} we showed that if $X$ is an
$n$-dimensional ALE manifold and $\beta\in(-n,-2)$,
then $\Delta:C^{k+2,\alpha}_{\beta+2}(X)\rightarrow 
C^{k,\alpha}_{\beta}(X)$ is an isomorphism. Our next
few results generalize this to a QALE manifold $X$.
We show that $\Delta:C^{k+2,\alpha}_{\beta,\gamma}(X)
\rightarrow C^{k,\alpha}_{\beta,\gamma-2}(X)$ is an 
isomorphism, for $(\beta,\gamma)$ in a certain nonempty 
open set ${\mathcal I}_X$ in $\R^2$. This is a useful 
analytic tool, and the principal justification for 
Definition~\ref{qholderdef}. 

We begin with a fairly weak existence result for solutions of 
the equation $\Delta u=f$. It can be proved by adapting known 
methods for ALE manifolds.

\begin{thm} Let\/ $(X,J,g)$ be a QALE K\"ahler manifold,  
$\alpha\in(0,1)$ and\/ $\beta,\gamma<0$. Then for each\/ 
$f\in C^{0,\alpha}_{\beta,\gamma-2}(X)$ there exists a 
unique $u\in C^{2,\alpha}(X)$ such that\/ $\Delta u=f$ 
and\/ $u(x)\rightarrow 0$ as $x\rightarrow\infty$ in~$X$.
\label{qlpexistthm}
\end{thm}

We quote a result about Schauder estimates on balls in $\R^n$,
taken from \cite[Th.~6.17]{GiTr}, that we will need next.

\begin{thm} Let\/ $B_1,B_2$ be the balls of radius $1,2$ in $\R^n$.
Let\/ $P$ be a linear elliptic operator of order 2 on functions on 
$B_2$, defined by
\begin{equation*}
Pu(x)=a^{ij}(x){\partial^2u\over\partial x^i\partial x^j}(x)
+b^i(x){\partial u\over\partial x^i}(x)+c(x)u(x).
\end{equation*}
Let\/ $k\ge 0$ be an integer and\/ $\alpha\in(0,1)$. Suppose the 
coefficients $a^{ij},b^i$ and\/ $c$ lie in $C^{k,\alpha}(B_2)$ 
and there are constants $\lambda,\Lambda>0$ such that 
$\bmd{a^{ij}(x)\xi_i\xi_j}\ge\lambda\ms{\xi}$ for all\/ 
$x\in B_2$ and $\xi\in\R^n$, and\/ $\cnm{a^{ij}}{k,\alpha}\le\Lambda$,
$\cnm{b^i}{k,\alpha}\le\Lambda$, and\/ $\cnm{c}{k,\alpha}\le\Lambda$
on $B_2$ for all\/ $i,j=1,\ldots,n$. Then there exists a constant\/ 
$C$ depending on $n,k,\alpha,\lambda$ and\/ $\Lambda$ such that 
whenever $u\in C^2(B_2)$ and\/ $f\in C^{k,\alpha}(B_2)$ with\/ 
$Pu=f$, we have $u\vert_{B_1}\in C^{k+2,\alpha}(B_1)$ and
\begin{equation*}
\bnm{u\vert_{B_1}}_{C^{k+2,\alpha}}\le 
C\bigl(\cnm{f}{k,\alpha}+\cnm{u}0\bigr).
\end{equation*}
\label{schauderthm}
\end{thm}

Using Theorems \ref{qlpexistthm} and \ref{schauderthm} we derive 
a sufficient condition for $\Delta:C^{k+2,\alpha}_{\beta,\gamma}(X)
\rightarrow C^{k,\alpha}_{\beta,\gamma-2}(X)$ to be an isomorphism.

\begin{thm} Let\/ $(X,J,g)$ be a QALE K\"ahler manifold, $(\rho,\sigma)$ 
a pair of radius functions on $X$, and\/ $\beta,\gamma<0$. Suppose 
there exists a smooth function $F:X\rightarrow(0,\infty)$ satisfying 
\begin{equation}
\Delta F\ge\rho^\beta\sigma^{\gamma-2}\quad\text{and}\quad
K_1\,\rho^\beta\sigma^\gamma\le F\le K_2\,\rho^\beta\sigma^\gamma
\label{qalefeq}
\end{equation}
for some $K_1,K_2>0$. Then whenever $k\ge 0$ and\/ 
$\alpha\in(0,1)$, there exists $C>0$ such that for each\/ 
$f\in C^{k,\alpha}_{\beta,\gamma-2}(X)$ there is a unique 
$u\in C^{k+2,\alpha}_{\beta,\gamma}(X)$ with\/ $\Delta u=f$, 
which satisfies $\nm{u}_{\smash{C^{k+2,\alpha}_{\beta,\gamma}}}\le 
C\nm{f}_{\smash{C^{k,\alpha}_{\beta,\gamma-2}}}$. In other words, 
$\Delta:C^{k+2,\alpha}_{\beta,\gamma}(X)\rightarrow 
C^{k,\alpha}_{\beta,\gamma-2}(X)$ is an isomorphism. 
\label{qlpisothm}
\end{thm}

\begin{proof} Let $f\in C^{k,\alpha}_{\beta,\gamma-2}(X)$, and suppose 
for simplicity that $\nm{f}_{\smash{C^{k,\alpha}_{\beta,\gamma-2}}}\le 1$. 
Then $f\in C^{0,\alpha}_{\beta,\gamma-2}(X)$, so by Theorem 
\ref{qlpexistthm} there exists a unique $u\in C^{2,\alpha}(X)$ such 
that $\Delta u=f$ and $u(x)\rightarrow 0$ as $x\rightarrow\infty$ in 
$X$. Since $\Delta u=f\le\rho^\beta\sigma^{\gamma-2}\le\Delta F$, we 
see that $\Delta(u-F)\le 0$ on $X$. Suppose that $u-F>0$ at some 
point of $X$. Then $u-F$ is non-constant and has a maximum in $X$, 
since $F>0$ and $u(x)\rightarrow 0$ as $x\rightarrow\infty$ in $X$. 
But this contradicts the maximum principle, as $\Delta(u-F)\le 0$. 
Therefore $u-F\le 0$, and $u\le F$. Similarly we show that $u\ge -F$, 
and so~$\md{u}\le F\le K_2\rho^\beta\sigma^\gamma$.

To complete the proof, it is enough to show that 
$u\in C^{k+2,\alpha}_{\beta,\gamma}(X)$ and 
$\nm{u}_{\smash{C^{k+2,\alpha}_{\beta,\gamma}}}\le C$ for 
some $C>0$ independent of $f$. We shall do this by applying 
Theorem \ref{schauderthm} to balls of radius $L\sigma(x)$ 
about $x$, for small~$L\in (0,{1\over 2})$. 

Let $B_1,B_2$ be the balls of radius 1 and 2 about 0 in 
$\R^{2m}$, where $m=\dim X$. Fix $x\in X$, and choose a vector 
space isometry $T_xX\cong\R^{2m}$. Let $L\in(0,{1\over 2})$ be 
small, and define a map $\Psi_x:B_2\rightarrow B_{2L\sigma(x)}(x)$
by $\Psi_x(y)=\exp_x\bigl(L\sigma(x)y\bigr)$, where 
$\exp_x:T_xX\rightarrow X$ is the exponential map, and we have 
identified $\R^{2m}$ and~$T_xX$.

Using the definition of QALE metric, we can show that
if $L>0$ is small, then
\begin{itemize}
\item $\Psi_x$ is a diffeomorphism between $B_2$ 
and $B_{2L\sigma(x)}(x)$,
\item $\Psi_x\vert_{\smash{B_1}}$ is a diffeomorphism 
between $B_1$ and $B_{L\sigma(x)}(x)$, and
\item the metric $g_x=L^{-2}\sigma(x)^{-2}\Psi_x^*(g)$ on $B_2$ 
is close to the Euclidean metric $h_0$ on $B_2$ in $C^{k+2,\alpha}$,
and $\nm{g_x-h_0}_{\smash{C^{k+2,\alpha}}}$ is bounded independently
of~$x\in X$.
\end{itemize}
The idea here is that $\sigma(x)$ is the approximate length-scale 
at which the metric on $X$ near $x$ differs from a Euclidean 
metric. Thus, balls of radius $L\sigma(x)$ and $2L\sigma(x)$ should
resemble Euclidean balls of the same radius provided $L$ is
sufficiently small.

Define an operator $P$ and functions $u',f'$ on $B_2$ by
\begin{equation*}
P=\bigl(L\sigma(x)\bigr)^2\Psi_x^*(\Delta), \quad
u'=\bigl(L\sigma(x)\bigr)^{-2}\Psi_x^*(u)
\quad\text{and}\quad f'=\Psi_x^*(f).
\end{equation*}
Then $P$ is the Laplacian of $g_x$ on $B_2$, which is close 
to $h_0$ in $C^{k+2,\alpha}$, and therefore $P$ is close in 
$C^{k,\alpha}$ to the Euclidean Laplacian $\Delta_0$ on $B_2$. 
Also $Pu'=f'$, as~$\Delta u=f$.

On $B_2$ we use the Euclidean metric $h_0$, and
on $B_{2L\sigma(x)}(x)$ we use the metric $g$. Since
$h_0\approx L^{-2}\sigma(x)^{-2}\Psi_x^*(g)$, it 
follows that
\begin{equation*}
\md{\nabla^jf'}\approx\bigl(L\sigma(x)\bigr)^j
\Psi_x^*\bigl(\md{\nabla^jf}\bigr)
\quad\text{and}\quad
\md{\nabla^ju'}\approx\bigl(L\sigma(x)\bigr)^{j-2}
\Psi_x^*\bigl(\md{\nabla^ju}\bigr).
\end{equation*}
Thus, as $\nabla^jf=O\bigl(\rho(x)^\beta\sigma(x)^{\gamma-2-j}\bigr)$ 
on $B_{2L\sigma(x)}$ we have $\nabla^jf'=
O\bigl(\rho(x)^\beta\sigma(x)^{\gamma-2}\bigr)$ on $B_2$, and
as $u=O\bigl(\rho(x)^\beta\sigma(x)^\gamma\bigr)$ on $B_{2L\sigma(x)}$ 
we have $u'=O\bigl(\rho(x)^\beta\sigma(x)^{\gamma-2}\bigr)$ 
on~$B_{2L\sigma(x)}$. 

Therefore Theorem \ref{schauderthm} shows that 
$\nabla^ju'=O\bigl(\rho(x)^\beta\sigma(x)^{\gamma-2}\bigr)$ on $B_1$
for $0\le j\le k+2$, together with the appropriate H\"older 
estimate. So $\nabla^ju=O\bigl(\rho(x)^\beta\sigma(x)^{\gamma-j}\bigr)$ 
on $B_{L\sigma(x)}(x)$ for $0\le j\le k+2$, together with the 
appropriate H\"older estimate. We have shown that $u$ is bounded
in $C^{k+2,\alpha}_{\beta,\gamma}(X)$, which completes the proof.
\end{proof}

Here when we say that $\Delta$ is an isomorphism between two 
Banach spaces of functions on $X$, we mean that it is an
{\it isomorphism of topological vector spaces}. That is, 
$\Delta$ is an isomorphism of vector spaces and is bounded 
and has bounded inverse, but it does not necessarily
identify the norms on the two Banach spaces. In fact 
\eq{qalefeq} is also necessary for 
$\Delta:C^{k+2,\alpha}_{\beta,\gamma}(X)\rightarrow 
C^{k,\alpha}_{\beta,\gamma-2}(X)$ to be an isomorphism.

\begin{prop} Let\/ $(X,J,g)$ be a QALE K\"ahler manifold, let\/ 
$(\rho,\sigma)$ be a pair of radius functions on $X$, and let\/ 
$\beta,\gamma<0$. Suppose that\/ $\Delta:C^{2,\alpha}_{\beta,\gamma}(X)
\rightarrow C^{0,\alpha}_{\beta,\gamma-2}(X)$ is an isomorphism. 
Then there exists a smooth function $F:X\rightarrow(0,\infty)$ 
satisfying~\eq{qalefeq}.
\label{qlpisoprop}
\end{prop}

\begin{proof} From \eq{rhoesteq} and \eq{siesteq} we see that 
$\rho^\beta\sigma^{\gamma-2}$ lies in $C^1_{\beta,\gamma-2}(X)$, 
and hence in $C^{0,\alpha}_{\beta,\gamma-2}(X)$. So there 
exists a unique $F\in C^{2,\alpha}_{\beta,\gamma}(X)$ with 
$\Delta F=\rho^\beta\sigma^{\gamma-2}$. Clearly $F$ is smooth, 
$\Delta F\ge\rho^\beta\sigma^{\gamma-2}$ and 
$F\le K_2\,\rho^\beta\sigma^\gamma$ for some $K_2>0$. 
Also, using \eq{rhoesteq} and \eq{siesteq} we can show that
$\Delta(K_1\rho^\beta\sigma^\gamma)\le\rho^\beta\sigma^{\gamma-2}$ 
for some small $K_1>0$, and thus by the proof of the previous 
theorem we have $K_1\rho^\beta\sigma^\gamma\le F$. Thus $F$ 
satisfies \eq{qalefeq}, as we want.
\end{proof}

\begin{dfn} Let $(X,J,g)$ be a QALE K\"ahler manifold, and define 
${\mathcal I}_X$ to be the set of pairs $(\beta,\gamma)\in\R^2$ such 
that $\beta<0$, $\gamma<0$ and $\Delta:C^{k+2,\alpha}_{\beta,\gamma}(X)
\rightarrow C^{k,\alpha}_{\beta,\gamma-2}(X)$ is an isomorphism for 
$k\ge 0$ and $\alpha\in(0,1)$. Theorem \ref{qlpisothm} and Proposition 
\ref{qlpisoprop} prove that this condition is independent 
of $k,\alpha$, and that $(\beta,\gamma)\in{\mathcal I}_X$ if and 
only if there exists a smooth function $F$ on $X$ satisfying 
\eq{qalefeq}. We can also show that ${\mathcal I}_X$ is an 
open set in~$\R^2$. 
\label{qixdef}
\end{dfn}

Of course these ideas are of no use at all if ${\mathcal I}_X$ is 
empty. In the following three results we show that it is not. 
The proof of the next proposition is elementary, and we omit it.

\begin{prop} Let\/ $\C^{m-n}$ be a subspace of\/ $\C^m$ for $n>0$,
and define $r,s:\C^m\rightarrow[0,\infty)$ by $r(x)=d(x,0)$ and\/ 
$s(x)=d(x,\C^{m-n})$. Let\/ $\beta,\gamma\in\R$. Then on 
$\C^m\setminus\C^{m-n}$ we have
\begin{equation}
\Delta\bigl(r^\beta s^\gamma)=
-{\textstyle{1\over 2}}r^{\beta-2}s^{\gamma-2}
\bigl[\gamma(\gamma+2n-2)r^2+\beta(2m-2+\beta+2\gamma)s^2\bigr].
\label{dersbegaeq}
\end{equation}
Since $0\le s\le r$, there exists $C>0$ such that\/ 
$\Delta\bigl(C\,r^\beta s^\gamma)\ge r^\beta s^{\gamma-2}$ 
if and only if
\begin{equation}
\gamma(\gamma+2n-2)<0\quad\text{and}\quad
\gamma(\gamma+2n-2)+\beta(2m-2+\beta+2\gamma)<0.
\label{begaineq}
\end{equation}
\label{begaineqprop}
\end{prop}

The pair of inequalities \eq{begaineq} are equivalent to
$2-2n<\gamma<0$ and
\begin{equation*}
\md{\beta+\gamma+m-1}<\sqrt{(m-1)^2+2\gamma(m-n)}.
\end{equation*}
For these to have a solution we must have $n>1$. Note also that
\begin{equation*}
(m-1)^2+2\gamma(m-n)=(m+1-2n)^2+2(\gamma+2n-2)(m-n). 
\end{equation*}
Thus if $2-2n<\gamma<0$ the square root 
$\sqrt{(m-1)^2+2\gamma(m-n)}$ exists, and 
\begin{equation*}
\md{m+1-2n}\le\sqrt{(m-1)^2+2\gamma(m-n)}\le m-1.
\end{equation*}
Hence, any solutions $\beta,\gamma$ to \eq{begaineq} must
satisfy $2-2m<\beta+\gamma<0$. Also, if $2-2n<\gamma<0$ and
$\beta+\gamma$ lies between $2-2n$ and $2(n-m)$ then 
$\beta,\gamma$ satisfy \eq{begaineq}. To get the factor 
${1\over 2}$ in \eq{dersbegaeq}, remember that $\Delta$ on 
K\"ahler manifolds is by convention half that on Riemannian 
manifolds.

Motivated by Proposition \ref{begaineqprop}, we can prove:

\begin{thm} Let\/ $(X,J,g)$ be a QALE K\"ahler manifold 
asymptotic to $\C^m/G$, and let\/ $n$ be the complex 
codimension of the singular set of\/ $\C^m/G$. Suppose 
$\beta,\gamma\in\R$ satisfy \eq{begaineq}. Then there exists a 
smooth function $F:X\rightarrow(0,\infty)$ satisfying~\eq{qalefeq}. 
\label{qfexistthm}
\end{thm}

The basic idea of the proof is to model $F$ on $\rho^\beta\sigma^\gamma$.
The details are complicated, and we will not give them. But here 
is a sketch of the case when $X$ is an ALE manifold. For 
$1-m<\delta<0$ we want to find a function $F$ such that 
$\Delta F\ge\rho^{2\delta-2}$ and $F\le C\,\rho^{2\delta}$. On 
$\C^m/G$ we have $\Delta(r^2)=-2m$ and $\bms{\nabla(r^2)}=4r^2$. 
We show that there exists a unique smooth function $u$ on $X$ 
such that $\Delta u=-2m$ and $u=\rho^2+O(\rho^{2-2m})$ for large 
$\rho$, and $4u-\ms{\nabla u}$ is bounded on~$X$.

Choose $K\in\R$ such that $u+K\ge 1$ and 
$\ms{\nabla u}\le 4(u+K)$ on $X$. Then
\begin{equation*}
\begin{split}
\Delta\bigl[(u+K)^\delta\bigr]
&=\delta(u+K)^{\delta-1}\Delta u
-{\textstyle{1\over 2}}\delta(\delta-1)(u+K)^{\delta-2}\ms{\nabla u}\\
&=-2\delta(\delta+m-1)(u+K)^{\delta-1}\\
&\qquad+{\textstyle{1\over 2}}\delta(\delta-1)
(u+K)^{\delta-2}\bigl(4(u+K)-\ms{\nabla u}\bigr)\\
&\ge -2\delta(\delta+m-1)(u+K)^{\delta-1},
\end{split}
\end{equation*}
since $\delta(\delta-1)>0$ and $\ms{\nabla u}\le 4(u+K)$. 
But there exist $C_1,C_2>0$ such that 
\begin{equation}
-2\delta(\delta+m-1)C_1(u+K)^{\delta-1}\ge\rho^{2\delta-2}
\quad\text{and}\quad
(u+K)^\delta\le C_2\rho^{2\delta}.
\end{equation}
So putting $F=C_1(u+K)^\delta$ we see that $\Delta F\ge\rho^{2\delta-2}$ 
and $F\le C_1C_2\rho^{2\delta}$, and the proof is finished. To 
generalize this proof to the case that $X$ is a QALE manifold,
we use similar functions $u$ on $X$ and on each $Y_i$ in the
decomposition of~$X$.

\begin{cor} Let\/ $(X,J,g)$ be a QALE K\"ahler manifold asymptotic 
to $\C^m/G$, and let\/ $n$ be the complex codimension of the 
singular set of\/ $\C^m/G$. Suppose that\/ $\beta,\gamma\in\R$ 
satisfy
\begin{equation}
\beta<0, \quad 2-2n<\gamma<0, \quad
\md{\beta+\gamma+m-1}<\sqrt{(m-1)^2+2\gamma(m-n)}.
\label{bgixeq}
\end{equation}
Then $(\beta,\gamma)\in{\mathcal I}_X$. Hence if\/ $2\le n\le m$ 
then ${\mathcal I}_X$ is nonempty.
\label{qixcor}
\end{cor}

To prove the last part, observe that if $2\le n\le m$ then 
\eq{bgixeq} holds for $\beta=1-m$ and $\gamma<0$ small. Note 
that $n\ge 2$ by definition, since $g$ is a QALE metric. If we 
allowed $n=1$ then the equation $2-2n<\gamma<0$ in \eq{bgixeq} 
would have no solutions, which illustrates the problems when 
$S$ has codimension one. 

By definition, if $(\beta,\gamma)\in{\mathcal I}_X$ then 
$\beta<0$ and $\gamma<0$. However, \eq{begaineq} admits 
solutions with $\gamma<0$ but $\beta\ge 0$. This suggests 
that $\Delta:C^{k+2,\alpha}_{\beta,\gamma}(X)\rightarrow 
C^{k,\alpha}_{\beta,\gamma-2}(X)$ can be an isomorphism 
when $\beta\ge 0$ as well. In fact this is the case, but we 
will not explore the idea because we shall need $\beta<0$ 
in our applications anyway. 

Using similar methods to those in the proof of Theorem 
\ref{qfexistthm}, we can show:

\begin{thm} Let\/ $(X,J,g)$ be a QALE K\"ahler manifold 
asymptotic to $\C^m/G$, and let\/ $\gamma<0$. Then there 
exists a smooth function $F:X\rightarrow(0,\infty)$ satisfying 
$\Delta F\ge\rho^{2-2m}\sigma^{\gamma-2}$ and\/ 
$K_1\rho^{2-2m}\le F\le K_2\rho^{2-2m}$ for some~$K_1,K_2>0$.
\label{q2-2mthm}
\end{thm}

From Theorem \ref{q2-2mthm} and the proof of Theorem 
\ref{qlpisothm} we deduce

\begin{cor} Let\/ $(X,J,g)$ be a QALE K\"ahler manifold of 
dimension $m$, and let $\gamma<0$. Then whenever $k\ge 0$ and\/ 
$\alpha\in(0,1)$, there exists $C>0$ such that for each\/ 
$f\in C^{k,\alpha}_{2-2m,\gamma-2}(X)$ there is a unique 
$u\in C^{k+2,\alpha}_{\smash{2-2m,0}}(X)$ with\/ $\Delta u=f$, 
which satisfies~$\nm{u}_{\smash{C^{k+2,\alpha}_{2-2m,0}}}\le 
C\nm{f}_{\smash{C^{k,\alpha}_{2-2m,\gamma-2}}}$. 
\label{q2-2mcor}
\end{cor}

The author believes that the following results are true.

\begin{conj} Let\/ $(X,J,g)$ be a QALE K\"ahler manifold 
asymptotic to $\C^m/G$, and let\/ $n$ be the complex codimension 
of the singular set of\/ $\C^m/G$. Then
\begin{equation}
{\mathcal I}_X=\bigl\{(\beta,\gamma)\in\R^2:\beta<0,\quad
2-2n<\gamma<0,\quad\beta+\gamma>2-2m\bigr\}.
\label{ixconjeq}
\end{equation}
Also, for generic $\beta,\gamma\in\R$ the map
$\Delta:C^{k+2,\alpha}_{\beta,\gamma}(X)\rightarrow 
C^{k,\alpha}_{\beta,\gamma-2}(X)$ is Fredholm, 
with finite-dimensional kernel and cokernel.
\end{conj}

\section{The Calabi conjecture for QALE manifolds}

The {\it Calabi conjecture} \cite{Cala} describes the possible
Ricci curvatures of K\"ahler metrics on a fixed compact complex
manifold $M$, in terms of the first Chern class $c_1(M)$. It was
proved by Yau \cite{Yau} in 1976. Since then several authors
such as Tian and Yau \cite{TiYa1,TiYa2} and Bando and Kobayashi
\cite{BaKo1,BaKo2} have proved versions of the Calabi conjecture
for {\it noncompact} manifolds.

We now state two versions of the Calabi conjecture for QALE 
manifolds, both of which will be proved in the author's book 
\cite{Joyc4}. So far as the author knows they do not follow
from any published noncompact version of the Calabi conjecture,
though there will obviously be similarities. The first, simpler 
version is based on the Calabi conjecture for ALE manifolds 
given in \cite[\S 6]{Joyc5}, and is proved in~\cite[\S 9.6]{Joyc4}.
\medskip

\noindent{\bf The Calabi conjecture for QALE manifolds (first version)} 
{\it Let\/ $(X,J,g)$ be a QALE K\"ahler manifold of dimension $m$,
with K\"ahler form $\omega$. Then
\begin{itemize}
\item[{\rm(a)}] Suppose that\/ $(\beta,\gamma)\in{\mathcal I}_X$, 
as in Definition \ref{qixdef}. Then for each\/ $f\in 
C^\infty_{\beta,\gamma-2}(X)$ there is a unique $\phi\in 
C^\infty_{\beta,\gamma}(X)$ such that\/ $\omega+\d\d^c\phi$ is a positive 
$(1,1)$-form and\/ $(\omega+\d\d^c\phi)^m={\rm e}^f\omega^m$ on~$X$.
\item[{\rm(b)}] Suppose $\gamma<0$. Then for each\/ $f\in 
C^\infty_{2-2m,\gamma-2}(X)$ there is a unique $\phi\in 
C^\infty_{2-2m,0}(X)$ such that\/ $\omega+\d\d^c\phi$ is a positive 
$(1,1)$-form and\/ $(\omega+\d\d^c\phi)^m={\rm e}^f\omega^m$ on~$X$.
\end{itemize}}
\medskip

We prove the conjecture by combining the method of Yau's
proof in the compact case \cite{Yau} with the ideas in \S 5
on analysis on QALE manifolds. The proof is similar to that 
of the Calabi conjecture for ALE manifolds given in 
\cite[\S 8.6-8.7]{Joyc4}, and discussed in~\cite[\S 6]{Joyc5}.

Here is a way to understand parts (a) and (b). Now 
$(\omega+\d\d^c\phi)^m={\rm e}^f\omega^m$ is really a nonlinear 
version of $\Delta\phi=-f$. Near infinity in $X$, where $\phi$ and 
$f$ are small, the two equations become very close. Therefore, if 
we can solve the equation $\Delta\phi=-f$ uniquely for $\phi,f$ in 
some given Banach spaces, we would na\"\i vely expect to be able to 
solve the Calabi conjecture in the same Banach spaces.

By Definition \ref{qixdef}, if $(\beta,\gamma)\in{\mathcal I}_X$ 
then $\Delta:C^{k+2,\alpha}_{\beta,\gamma}(X)\rightarrow 
C^{k,\alpha}_{\beta,\gamma-2}(X)$ is an isomorphism. 
Thus we can solve the equation $\Delta\phi=-f$ uniquely 
with $\phi\in C^{k+2,\alpha}_{\beta,\gamma}(X)$ and
$f\in C^{k,\alpha}_{\beta,\gamma-2}(X)$, and 
similarly with $\phi\in C^\infty_{\beta,\gamma}(X)$ 
and~$f\in C^\infty_{\beta,\gamma-2}(X)$. 

So part (a) above, which says we can solve the Calabi 
conjecture with $\phi\in C^\infty_{\beta,\gamma}(X)$ and
$f\in C^\infty_{\beta,\gamma-2}(X)$, is as we would expect.
Similarly, part (b) above corresponds to Corollary 
\ref{q2-2mcor}, which says that we can solve $\Delta\phi=-f$ 
uniquely with $\phi\in C^{k+2,\alpha}_{2-2m,0}(X)$ and
$f\in C^{k,\alpha}_{2-2m,\gamma-2}(X)$, and so with $\phi\in 
C^\infty_{2-2m,0}(X)$ and~$f\in C^\infty_{2-2m,\gamma-2}(X)$.

So far we have estimated the decay of functions $\phi$ and 
metrics $g$ on $X$ in two different ways. In \S 2-\S 4 we 
pulled $g$ and $\phi$ back to $V_i\times Y_i$ and bounded 
them using powers of functions $\mu_{i,j}$ and $\nu_i$ on 
$V_i\times Y_i$. In \S 5, and in the conjecture above and 
its proof, we estimated $\phi$ in terms of powers of 
functions $\rho,\sigma$ on~$X$. 

There are good reasons for using these two approaches in the 
way we have --- attempting to define QALE metrics directly on 
$X$, and trying to solve the Calabi conjecture by pulling back 
to $V_i\times Y_i$, both seem to lead to disaster. However, in 
order to construct QALE K\"ahler metrics with prescribed Ricci 
curvature on $X$ we need to integrate these two approaches, 
because in solving the Calabi conjecture on $X$ we need $\phi$
to be a function of {\it K\"ahler potential type}, as in~\S 4.

Here is a version of the Calabi conjecture for QALE manifolds 
which achieves this, which is proved in~\cite[\S 9.7]{Joyc4}.
\medskip

\noindent{\bf The Calabi conjecture for QALE manifolds (second version)} 
{\it Let\/ $(X,J,g)$ be a QALE K\"ahler manifold asymptotic to $\C^m/G$
with\/ $\Fix(G)=\{0\}$, let\/ $\omega$ be the K\"ahler form and\/ $\xi$ 
the Ricci form of $g$, and let\/ $\epsilon<-2$. Suppose $f$ is a smooth 
function on $X$ with
\begin{equation}
\nabla^l\psi_i^*(f)=\sum\begin{Sb}j\in I:i\nsucceq j\end{Sb}
O\bigl(\mu_{i,j}^{d_j}\nu_i^{\epsilon-l}\bigr)
\label{qfmneq}
\end{equation}
on $V_i\times Y_i\,\big\backslash U_i$, for all\/ 
$i\in I\setminus\{\infty\}$ and\/ $l\ge 0$. Then there is a unique 
smooth function $\phi$ on $X$ such that\/ $\omega'=\omega+\d\d^c\phi$ 
is a positive $(1,1)$-form and\/ $(\omega')^m=e^f\omega^m$ on $X$, and 
\begin{equation}
\nabla^l\psi_i^*(\phi)=\sum\begin{Sb}j\in I:i\nsucceq j\end{Sb}
O\bigl(\mu_{i,j}^{d_j}\nu_i^{-l}\bigr)
\label{qphimneq}
\end{equation}
on $V_i\times Y_i\,\big\backslash U_i$, for all\/ 
$i\in I\setminus\{\infty\}$ and\/ $l\ge 0$. Furthermore, $\phi$ is 
of K\"ahler potential type on $X$, and the K\"ahler metric $g'$ on 
$X$ with K\"ahler form $\omega'$ is QALE and has Ricci 
form~$\xi'=\xi-{1\over 2}\d\d^c f$.}
\medskip

We prove this by translating \eq{qfmneq} and \eq{qphimneq}
into equations directly on $X$, prescribing the asymptotic
behaviour of $\phi$ and $f$ in terms of $\rho,\sigma$ and
some similar functions $\rho_i$ for $i\in I\setminus\{0\}$.
Then we apply the method used to prove the first version
above --- this is really a more complicated version of
part (b) of the first conjecture.

The point of assuming that $\Fix(G)=V_\infty=\{0\}$ here is 
that taken over all $i\in I\setminus\{\infty\}$, eqns \eq{qfmneq} 
and \eq{qphimneq} prescribe the asymptotic behaviour of 
$f$ and $\phi$ everywhere on $X$, except within a fixed
distance of $\pi^{-1}(V_\infty)$. If $\Fix(G)=\{0\}$ this 
is all of $X$ except a compact subset. But if 
$\Fix(G)\ne\{0\}$ then $\pi^{-1}(V_\infty)$ extends to 
infinity in $X$, and \eq{qfmneq} and \eq{qphimneq} only 
prescribe the behaviour on a part of the `boundary' of 
$X$, so they are not good boundary conditions.

\section{The proofs of Theorems \ref{qalerfthm} and \ref{qaleholthm}}

Before proving the theorems, we give two preliminary propositions.

\begin{prop} Suppose $(X,J,g)$ is a Ricci-flat QALE K\"ahler manifold. 
Then $g$ is the only Ricci-flat QALE K\"ahler metric in its K\"ahler class.
\label{qrf1prop}
\end{prop}

\begin{proof} Suppose for a contradiction that $g,g'$ are 
distinct Ricci-flat QALE K\"ahler metrics on $X$ in the same K\"ahler 
class, and let $X$ be of the smallest dimension in which this 
can happen. Clearly $\Fix(G)=\{0\}$, since otherwise we can
replace $X$ by $Y_\infty$, which has smaller dimension. Let 
$g,g'$ be asymptotic to metrics $g_i,g_i'$ on $Y_i$ for 
$i\in I$, as in Definition \ref{qaledef}. Then $g_i,g_i'$ 
are Ricci-flat, and in the same K\"ahler class. Thus $g_i=g_i'$ for 
all $i\ne\infty$ in $I$, since $\dim Y_i<\dim X$ when~$i\ne\infty$.

Let $\omega,\omega'$ be the K\"ahler forms of $g,g'$. Then by Theorem 
\ref{kpotqalethm} we have $\omega'=\omega+\d\d^c\phi$, where $\phi$ 
is a function of K\"ahler potential type on $X$. Since $g_i=g_i'$ 
for $i\ne\infty$, the functions $\phi_i$ of Definition 
\ref{qktypedef} are zero for $i\ne\infty$, and $\phi$ satisfies 
\eq{qphimneq}. But $g,g'$ have the same Ricci curvature, so 
that $(\omega')^m=\omega^m$. Thus $\phi$ and $f\equiv 0$ satisfy 
the second version of the Calabi conjecture for QALE manifolds.
By uniqueness in the conjecture we see that $\phi\equiv 0$, 
so that $g=g'$, a contradiction. Thus $g$ is unique in its 
K\"ahler class.
\end{proof}

\begin{prop} Let\/ $(X,\pi)$ be a crepant resolution of\/ 
$\C^m/G$ with\/ $\Fix(G)=\{0\}$, and let\/ $g$ be a QALE K\"ahler 
metric on $X$ asymptotic to metrics $g_i$ on $Y_i$ for $i\in I$, 
such that\/ $\psi_i^*(g)=h_i\times g_i$ on $V_i\times Y_i\setminus U_i$. 
Suppose that for each\/ $i\in I\setminus\{\infty\}$ there is a Ricci-flat 
K\"ahler metric $g_i'$ on $Y_i$ in the same K\"ahler class as $g_i$. Then 
there exists a QALE K\"ahler metric $g'$ on $X$ in the same K\"ahler class 
as $g$ with Ricci form ${1\over 2}\d\d^cf$, where $f$ is a smooth 
function on $X$ satisfying
\begin{equation}
\nabla^l\psi_i^*(f)=
\sum\begin{Sb}j\ne k\in I\setminus\{\infty\}:\\ V_j\cap V_k=\{0\}\end{Sb}
O\bigl(\mu_{i,j}^{\smash{d_j}}\mu_{i,k}^{\smash{d_k}}\nu_i^{-4-l}\bigr)
\quad\text{on $V_i\times Y_i\,\big\backslash U_i$,}
\label{nabfmneq}
\end{equation}
for $i\in I$ and\/ $l\ge 0$. However, if there are no $j,k\in 
I\setminus\{\infty\}$ with\/ $V_j\cap V_k=\{0\}$ then \eq{nabfmneq} 
does not hold, but instead\/ $f$ is compactly-supported on~$X$.
\label{qrf2prop}
\end{prop}

\begin{proof} Let $g_i,g_i'$ have K\"ahler forms $\omega_i,\omega_i'$. Then 
Theorem \ref{kpotqalethm} shows that $\omega_i'=\omega_i+\d\d^c\phi_i$, 
where $\phi_i$ is a unique function of K\"ahler potential type on $Y_i$. 
As $g_i'$ is the only Ricci-flat QALE metric in its K\"ahler class by 
Proposition \ref{qrf1prop}, we see that $\chi_{\gamma,i}^*
(g_{\gamma\cdot i}')=g_i'$, and $g_i'$ is asymptotic to $g_j'$ 
when $i,j\ne\infty$ and $i\succeq j$. Using these facts one can 
show that the $\phi_i$ satisfy parts (i) and (ii) of 
Proposition~\ref{qkpotcondprop}. 

Therefore we may apply Theorem \ref{phiextthm} and Proposition 
\ref{phiextprop} to find a function $\phi$ of K\"ahler potential type 
on $X$ such that $\omega'=\omega+\d\d^c\phi$ is the K\"ahler form of a 
QALE K\"ahler metric $g'$ on $X$, which is asymptotic to the Ricci-flat
metrics $g_i'$ on $Y_i$. Also, as $\Fix(G)=\{0\}$, by construction 
$\phi$ satisfies \eq{qalephiieq} outside a compact subset $T$ of~$X$. 

Let $\omega_0$ be the K\"ahler form of the Euclidean metric $h_0$ 
on $\C^m/G$, and let $\Omega_0=\d z^1\wedge\cdots\wedge\d z^m$ be 
the holomorphic volume form on $\C^m/G$, which is well-defined 
as $G\subset{\rm SU}(m)$, since $X$ is a crepant resolution. 
Calculation shows that $\omega_0^m=C_m\Omega_0\wedge\overline\Omega_0$ 
on $\C^m/G$, where $C_m=2^{-m}i^mm!(-1)^{m(m-1)/2}$. Let 
$\Omega=\pi^*(\Omega_0)$. Then $\Omega$ is a nonsingular holomorphic 
volume form on $X$. Define a smooth real function $f$ on $X$ 
by ${\rm e}^f(\omega')^m=C_m\Omega\wedge\overline\Omega$. So $g'$ 
has Ricci form ${1\over 2}\d\d^cf$, as we want.

We must show that $f$ satisfies \eq{nabfmneq}. For simplicity 
we will restrict our attention to the case $N(V_i)=G$ for all 
$i\in I$. Then \eq{qalephijeq} gives
\begin{equation}
\phi=\sum_{i\in I\setminus\{\infty\}}k_i\Phi_i'
\quad\text{on $X\setminus T$.}
\label{phiphiieq}
\end{equation}
Let $x\in X$ satisfy $d\bigl(\pi(x),S\bigr)>2R$. Then putting 
$i=0$ in $\psi_i^*(g)=h_i\times g_i$ on $V_i\times Y_i\setminus U_i$ 
shows that $g=\pi^*(h_0)$ near $x$. Thus $\omega=\pi^*(\omega_0)$ and 
$\Omega=\pi^*(\Omega)$, and so $\omega^m=C_m\Omega\wedge\overline\Omega$ 
near~$x$. 

If $i\ne\infty$ in $I$ and $(v,y)\in V_i\times Y_i\setminus U_i$ with 
$\psi_i(v,y)=x$, then $\psi_i^*(g)=h_i\times g_i$ near $(v,y)$. Also 
the function $\Phi_i$ defined in \eq{phiipdefeq} satisfies $\Phi_i=\phi_i$ 
near $(v,y)$, since $\mu_{i,j}>2R$ near $(v,y)$ for all $j\ne 0$. But 
$N(V_i)=G$, and so $\Phi_i=\psi_i^*(\Phi_i')$, where $\Phi_i'$ is 
defined by \eq{qalephijeq}. Therefore the metric $h_i\times g_i'$ on 
$V_i\times Y_i$ has K\"ahler form $\psi_i^*(\omega+\d\d^c\Phi_i')$ 
near~$(v,y)$. 

But $h_i\times g_i'$ is Ricci-flat. So $\omega+\d\d^c\Phi_i$ is the 
K\"ahler form of a Ricci-flat metric near $x$. It is then not difficult 
to show that $(\omega+\d\d^c\Phi_i')^m=C_m\Omega\wedge\overline\Omega$ 
near $x$. Thus we have ${\rm e}^f(\omega+\d\d^c\phi)^m=\omega^m=
(\omega+\d\d^c\Phi_i')^m$ wherever $d\bigl(\pi(x),S\bigr)>2R$ 
on $X$, for all $i\in I\setminus\{\infty\}$. Multiplying out 
$(\omega+\d\d^c\Phi_i')^m=\omega^m$ and rearranging gives
\begin{equation}
m\,\d\d^c\Phi'_i\wedge\omega^{m-1}=-{\textstyle{1\over 2}}
m(m-1)(\d\d^c\Phi_i')^2\wedge\omega^{m-2}+\ldots,
\label{mddcphieq}
\end{equation}
where `\dots' are terms of order at least 3 in $\d\d^c\Phi_i'$. 
Substituting \eq{phiphiieq} into the equation 
${\rm e}^f(\omega+\d\d^c\phi)^m=\omega^m$, rearranging and using 
\eq{mddcphieq} shows that
\begin{equation}
\begin{split}
{-2f\over m(m-1)}\,\omega^m=&\sum
\begin{Sb}i\ne j\in I\setminus\{\infty\}\end{Sb}
k_ik_j\,\d\d^c\Phi_i'\wedge\d\d^c\Phi_j'\wedge\omega^{m-2}\\
+&\sum\begin{Sb}i\in I\setminus\{\infty\}\end{Sb}
k_i(k_i-1)(\d\d^c\Phi_i')^2
\wedge\omega^{m-2}+\ldots,
\end{split}
\label{efommeq}
\end{equation}
where `\dots' are terms of order at least 3 in the $\d\d^c\Phi_i'$ 
or at least 2 in $f$, and the equation holds for $x\in X\setminus T$ 
with $d\bigl(\pi(x),S\bigr)>2R$. Thus $f$ is roughly quadratic 
in the $\d\d^c\Phi_i'$, to highest order. 

Using \eq{efommeq} and the special properties of the $k_i$ and 
$\Phi_i'$, we can show that \eq{nabfmneq} holds. The proof of
this is complicated, and we will not give it. The basic idea 
is that $g'$ is made by combining the Ricci-flat metrics 
$h_i\times g_i'$. The dominant terms in $\Ric(g')$ result from 
interference between $h_j\times g_j'$ and $h_k\times g_k'$ for 
$j\ne k$, and contribute $O(\mu_{i,j}^{\smash{d_j}}\nu_i^{-2}
\cdot\mu_{i,k}^{\smash{d_k}}\nu_i^{-2})$ to $\psi_i^*(f)$ on 
$V_i\times Y_i\setminus U_i$, and corresponding terms to 
$\nabla^l\psi_i^*(f)$.

However, if $V_j\cap V_k=V_i$ with $i\ne\infty$ then 
$h_i\times g_i'$ is Ricci-flat and asymptotic to both 
$h_j\times g_j'$ and $h_k\times g_k'$. So we introduce no 
extra Ricci curvature by combining $h_j\times g_j'$ and 
$h_k\times g_k'$ in this case, which is why the sum in 
\eq{nabfmneq} is restricted to $j,k$ with $V_j\cap V_k=\{0\}$.
More details of this argument are given in~\cite[\S 9.8]{Joyc4}.
\end{proof}

\subsection{Proof of Theorem \ref{qalerfthm}}

We work by induction on $m=\dim X$. The result is trivial for 
$m=0,1$, giving the first step. For the inductive step, suppose 
that $X$ is a crepant resolution of $\C^m/G$ for some $m\ge 2$, 
that $\kappa$ is a K\"ahler class on $X$ containing QALE K\"ahler 
metrics, and that the theorem is true in dimensions $0,1,\ldots,m-1$. 
We shall show that there exists a Ricci-flat QALE K\"ahler metric 
$\hat g$ on $X$ in $\kappa$. This $\hat g$ is unique by Proposition 
\ref{qrf1prop}, and so each K\"ahler class on $X$ contains a unique 
Ricci-flat QALE K\"ahler metric. Thus by induction the theorem is 
true for all~$m$.

If $\Fix(G)\ne\{0\}$ then $X=V_\infty\times Y_\infty$, and as 
$\dim Y_\infty<m$ there is a unique Ricci-flat QALE K\"ahler metric 
$\hat g_\infty$ on $Y_\infty$ in $\kappa\vert_{\smash{Y_\infty}}$.
Then $\hat g=h_\infty\times\hat g_\infty$ is the metric on $X$ that
we seek. So suppose that~$\Fix(G)=\{0\}$.

By Theorem \ref{qalekthm} we can choose QALE K\"ahler metrics 
$g$ on $X$ and $g_i$ on $Y_i$ such that $g$ has K\"ahler class 
$\kappa$ and $\psi_i^*(g)=h_i\times g_i$ on $V_i\times Y_i\setminus U_i$,
where $U_i$ is defined using some $R>0$ depending on $\kappa$.
Now $Y_i$ is a crepant resolution of $W_i/A_i$, and 
$\dim Y_i<m=\dim X$ if $i\ne\infty$. Thus for each $i\ne\infty$ 
in $I$ there is a unique Ricci-flat K\"ahler metric $g_i'$ on $Y_i$ 
in the K\"ahler class of $g_i$, by the inductive hypothesis.

Proposition \ref{qrf2prop} applies, and gives a QALE K\"ahler 
metric $g'$ on $X$ in the K\"ahler class $\kappa$ with Ricci form
${1\over 2}\d\d^cf$, where either $f$ satisfies \eq{nabfmneq} if 
there exist $j,k\in I\setminus\{\infty\}$ with $V_j\cap V_k=\{0\}$, 
or else $f$ is compactly-supported on $X$. Now it is easy
to show from \eq{nabfmneq} that $f$ satisfies \eq{qfmneq} 
with~$\epsilon=-4$. 

Therefore by the second version of the Calabi conjecture 
for QALE manifolds, which holds by \cite[\S 9.7]{Joyc4},
there is a unique function $\phi'$ on $X$ such that 
$\hat\omega=\omega'+\d\d^c\phi'$ satisfies $\hat\omega^m=
{\rm e}^f(\omega')^m$ and \eq{qphimneq} holds for $\phi'$. 
The K\"ahler metric $\hat g$ on $X$ with K\"ahler form 
$\hat\omega$ is QALE, as $\phi'$ is of K\"ahler potential 
type, and has Ricci form $\xi-{1\over 2}\d\d^cf$, where $\xi$ 
is the Ricci form of $g'$. But $\xi={1\over 2}\d\d^cf$ from
above, and so $\hat g$ is Ricci-flat. Thus we have found a 
Ricci-flat QALE K\"ahler metric $\hat g$ in the K\"ahler class 
$\kappa$, and the proof is complete.

\subsection{Proof of Theorem \ref{qaleholthm}}

Let $(X,J,g)$ be a Ricci-flat QALE K\"ahler manifold asymptotic 
to $\C^m/G$, where $\C^m/G$ is irreducible. As $X$ is 
simply-connected and $g$ is Ricci-flat K\"ahler we have 
$\Hol(g)\subseteq{\rm SU}(m)$. But $g$ is nonsymmetric as it is 
Ricci-flat, and since $\C^m/G$ is irreducible one can show 
that $g$ is irreducible. So Berger's classification of
Riemannian holonomy groups \cite[\S 10]{Sala} shows
that $\Hol(g)$ is ${\rm SU}(m)$ or ${\rm Sp}(m/2)$. The Euclidean 
metric $h_0$ on $\C^m/G$ has $\Hol(h_0)=G$, and as $g$ is 
asymptotic to $h_0$ we see that $G\subset\Hol(g)$. Thus, 
if $G$ is not conjugate to a subgroup of ${\rm Sp}(m/2)$ then 
$\Hol(g)\ne{\rm Sp}(m/2)$, which forces~$\Hol(g)={\rm SU}(m)$.

So suppose $m\ge 4$ is even and $G$ is conjugate to a subgroup 
of ${\rm Sp}(m/2)$. Then there exists a $G$-invariant, constant 
complex symplectic form $\omega_{\scriptscriptstyle\mathbb C}$ 
on $\C^m$, which pushes down to $\C^m/G$. The pull-back 
$\pi^*(\omega_{\scriptscriptstyle\mathbb C})$ is a nonsingular 
complex symplectic form on $X$, and using a Bochner argument we 
can prove that $\nabla\pi^*(\omega_{\scriptscriptstyle\mathbb C})=0$. 
Therefore $\Hol(g)\subset{\rm Sp}(m/2)$, and so~$\Hol(g)={\rm Sp}(m/2)$.

\end{document}